\documentclass[10pt,draft]{amsart}

\usepackage{amssymb,amsfonts}
\usepackage{euscript}

\usepackage[ps,matrix,arrow,curve,cmtip]{xy}
\newcommand{\margnote}[1]{\mbox{}\marginpar{\tiny\hspace{0pt}#1}}

\numberwithin{equation}{section}

\makeatletter
  \let\c@subsection\c@equation
\makeatother

\theoremstyle{plain}   %% This is the default, anyway

\newtheorem{thm}[subsection]{Theorem}
\newtheorem{prop}[subsection]{Proposition}
\newtheorem{cor}[subsection]{Corollary}
\newtheorem{lemma}[subsection]{Lemma}

\theoremstyle{remark}
\newtheorem{rem}[subsection]{Remark}    
\newtheorem{exam}[subsection]{Example}

\newtheorem{defn}[subsection]{Definition}

\theoremstyle{plain}

\DeclareMathOperator{\id}{id}
\DeclareMathOperator{\colim}{colim}

\newcommand{\op}{{\operatorname{op}}}
\newcommand{\ob}{{\operatorname{ob}}}

\newcommand{\ra}{\rightarrow}

\newcommand{\xra}{\xrightarrow}

\newcommand{\cat}[1]{{\operatorname{\mathbf{#1}}}}

\newcommand{\Set}{{\operatorname{\EuScript{S}et}}}
\newcommand{\sSet}{{\operatorname{\EuScript{S}}}}

\DeclareMathOperator{\Tor}{Tor}

\newcommand{\Z}{\mathbb{Z}}
\newcommand{\N}{\mathbb{N}}

\newcommand{\Q}{\mathbb{Q}}

\newcommand{\dfn}{\textbf}

\def\noloc{\;{:}\,}

\newcommand{\forcepar}{\mbox{}\par}

\title{Every homotopy theory of simplicial  algebras admits a proper
model}   
\author{Charles Rezk}
\date{March 10, 2000}
\address{Institute for Advanced Study \\
Princeton, NJ 08540}
\email{rezk@math.nwu.edu}
\thanks{The author was supported by an AMS Centennial Fellowship.}

\renewcommand{\margnote}[1]{}

\renewcommand{\Set}{{\EuScript{S}}}
\renewcommand{\sSet}{{s\Set}}
\renewcommand{\phi}{\varphi}

\newcommand{\indset}[1]{{\EuScript{#1}}}
\newcommand{\fsets}{{f\Set}}
\newcommand{\fsetsind}[1]{{f\Set(\indset{#1})}}
\newcommand{\fsetsdind}[2]{{f\Set(\indset{#1},\indset{#2})}}
\newcommand{\fsetsgr}[1]{{\N(\indset{#1})}}
\newcommand{\fset}[1]{{{#1}}}

\newcommand{\Endth}[1]{{{\mathcal{E}}_{#1}}}
\newcommand{\theories}{{\mathcal{T}}}
\newcommand{\alg}[1]{{#1\text{-}\mathrm{alg}}}
\newcommand{\bimod}[1]{{#1\text{-}\mathrm{mod}}}
\newcommand{\Ex}{{\operatorname{Ex}}}

\newcommand{\rtmod}[1]{\bimod{I,#1}}
\newcommand{\ev}{{\operatorname{ev}}}
\newcommand{\Id}{{\operatorname{Id}}}
\newcommand{\Sing}{{\operatorname{Sing}}}
\newcommand{\fgf}{{\operatorname{fgf}}}

\newcommand{\freetheory}{{\mathcal{F}}}
\newcommand{\freeser}{{\mathcal{S}}}

\newcommand{\Func}{{\operatorname{Func}}}
\newcommand{\filtfunc}{{\operatorname{Func}^{\operatorname{fr}}}}
\newcommand{\filtofunc}{{\operatorname{Func}^{\operatorname{f}}}}
\newcommand{\dfiltfunc}{{\operatorname{Func}^{\operatorname{dfr}}}}
\newcommand{\degencat}{\Delta_{+}}
\newcommand{\simpind}{{\Delta}}
\newcommand{\freetrees}{{\mathcal{Q}}}
\newcommand{\esstrees}{{\mathcal{Q}_e}}
\newcommand{\inputs}{{\operatorname{in}}}
\newcommand{\outputs}{{\operatorname{out}}}
\newcommand{\tree}[1]{\mathcal{#1}}

\newcommand{\Ho}{{\operatorname{Ho}}\,}
\newcommand{\diag}{{\operatorname{diag}}}

\theoremstyle{plain}

\newtheorem*{thm-a}{Theorem A}
\newtheorem*{thm-b}{Theorem B}
\newtheorem*{thm-c}{Theorem C}
\newtheorem*{thm-d}{Theorem D}
\newtheorem*{cor*}{Corollary}

\begin{document}

%%% abstract
\begin{abstract}
We show that any closed model category of simplicial algebras over an
algebraic theory is
Quillen equivalent to a proper closed model category.   
By ``simplicial algebra'' we mean any category of algebras over a
simplicial algebraic 
theory, which is allowed to be multi-sorted.  The results have
applications to the construction of localization model category
structures. 
\end{abstract}

%%% the title
\maketitle

%%% table of contents
%\tableofcontents

\section{Introduction}

To axiomatize the notion of a ``homotopy theory'' Quillen introduced
closed model categories 
\cite{homotopical-algebra}, and produced a number of examples of such, one
class of which are categories of \emph{simplicial algebras}.
A standard technique for constructing new model categories from old ones
is that of 
\emph{localization}: given a category 
$\cat{C}$ equipped with a model
category structure and a morphism $f$ in that category, one produces a
new model category structure on $\cat{C}$ in which the weak
equivalences are the smallest class containing both the old weak
equivalences and the map $f$.  There are several ``machines'' for
constructing localization model category structures; one of the most general
is due to Hirschhorn~\cite{hirschhorn}; note also
\cite{abstract-homotopy}, \cite{smith-combinatorial-model-cats}, and
\cite[Ch.\ X]{goerss-jardine-simplicial-book}.  They have been used
extensively in recent 
years, notably\margnote{cite spectrum people} to construct model
categories for stable homotopy theories. 

These localization machines require that the initial model category
structure on $\cat{C}$ 
have certain additional properties, beyond those introduced by
Quillen.  In most cases they require in particular that $\cat{C}$ be
a ``left proper'' model category; namely, the class of weak
equivalences should be closed under cobase change along cofibrations
(see \eqref{subsec-proper-def}).  Properness was first introduced by
Bousfield and Friedlander \cite{bousfield-friedlander} as an axiom
needed to be able to put a model structure on a category of spectra;
their construction is in fact an instance of a localization model
category.  

Many well-understood examples of model categories, including most
categories of simplicial algebras, turn out \emph{not} to be
proper.  (We give examples of such in
\S\ref{sec-improper-model-cat}.)  It is the most non-trivial 
axiom needed for the localization machines.  Thus, the following
question becomes significant: does our homotopy theory admit a \emph{proper}
model? 
That is, given a closed model category $\cat{C}$
which is not necessarily proper, does there exist a proper
closed model category $\cat{C}'$ which has the same homotopy theory as
$\cat{C}$?  

In this paper, we examine the case of \emph{simplicial
algebras}, i.e., simplicial objects in a category of algebras
associated to an algebraic theory in the sense of Lawvere
\cite{lawvere-functorial-semantics}, and more generally the case of simplicial
algebras over a \emph{multi-sorted, simplicial theory} 
(see \S\ref{subsec-theories}, \S\ref{subsec-simplicial-alg}).  This
class of examples includes simplicial groups, rings, and so forth, as
well as algebras over a simplicial operad, as in \cite{thesis}.  They
are the simplicial analogues of the topological theories considered by
Boardman and Vogt \cite{boardman-vogt-homotopy-invariant-structures}.
Categories of simplicial algebras always  always admit a 
model category structure, with the weak equivalences being those of
the underlying simplicial sets
\eqref{thm-closed-model-cat}.  
\begin{thm-a}\label{thm-a}
The homotopy theory of a category of simplicial algebras always admits
a proper model.
\end{thm-a}
Whether \emph{any} reasonable homotopy theory
(e.g., one associated to a model category) admits a proper model is an
open question;
Theorem A is the only result in this direction that I am aware of.

Theorem A can be made more precise.  It is a corollary of the
following 
\begin{thm-b}\label{thm-b}
Let $T$ be a (possibly simplicial, possibly multi-sorted) theory, and let
$\alg{T}$ be the corresponding 
category of simplicial $T$-algebras, equipped with a
simplicial model category
structure in which a 
map is a weak equivalence or fibration if it is a weak equivalence or
fibration of the
underlying simplicial sets.

Then there exists a morphism $S\ra T$ of simplicial theories such
that
\begin{enumerate}
\item [(1)] the induced adjoint pair $\alg{S}\rightleftarrows\alg{T}$
is a Quillen equivalence of model categories, and

\item [(2)] $\alg{S}$ is a \emph{proper} simplicial closed model category.
\end{enumerate}
\end{thm-b}
The proof of Theorem B follows a straightforward pattern; we (a) put a
model category structure on the category of simplicial theories
\eqref{thm-closed-model-cat} so that in particular cofibrant
resolutions of simplicial theories exist; (b) show 
that algebras over a \emph{cofibrant} 
simplicial theory are a proper model category
\eqref{cor-cof-theory-proper-model-cat}; and (c) observe that weakly
equivalent simplicial theories give rise to Quillen equivalent model
categories 
of algebras \eqref{cor-ho-inv-of-simpl-algebras}.

We say that a category is \emph{pointed} if the initial object is
isomorphic to the terminal object.  It is most natural to study stable
homotopy of algebras in the context of pointed objects.  Thus we offer 
\begin{thm-c}\label{thm-c}
Given the hypotheses of Theorem B, suppose that in addition
$\alg{T}$ is a pointed category.  Then $S$ can be chosen as in
Theorem B so that $\alg{S}$ is also a pointed
category. 
\end{thm-c}

Finally, we have
\begin{thm-d}
Given the hypotheses of Theorem B (resp.\ of Theorem C), the theory
$S$ can be chosen as in Theorem B (or Theorem C) so that $\alg{S}$ is
a \emph{cellular model category} in the sense of Hirschhorn
\cite{hirschhorn} .  
\end{thm-d}
By Hirschhorn's results \cite{hirschhorn}, Theorem D implies 
\begin{cor*}
For any set of
maps in $\alg{S}$ there is a localization model category structure
with respect to this set.
\end{cor*}
The proofs of Theorems A, B, C and D are given in
\S\ref{sec-proofs-of-theorems}. 

In order to prove these results, we need to set up a certain amount of
foundations for algebraic theories and their homotopy theory; this
will take all of
\S\S\ref{sec-functors-on-sets}--\ref{sec-ho-invariance-prop}.  Our
exposition of theories
(\S\S\ref{sec-functors-on-sets}--\ref{sec-theories-algebras-bimod})
is more involved than one might like; this is because we want to deal
with ``multi-sorted'' theories, and because we need to introduce the
notion of ``bimodules'' of algebraic theories.  However, this is not idle
generalization: the category of single-sorted theories and categories
of bimodules over such are themselves categories of algebras over a
multi-sorted theory, 
so considering multi-sorted theories from the
start lets us avoid much duplication of exposition.  The theory of
bimodules of 
algebraic theories plays
an important role in the proofs of the main theorems (see
\S\ref{sec-ho-invariance-prop} and \S\ref{sec-cofs-of-theories}).

Some of this foundational material seems to be of independent
interest, notably our definition of 
``bimodules'' of algebraic theories and their relation to functors
between categories of algebras \eqref{subsec-bimodules-def}, and
the homotopy invariance results
\eqref{rem-derived-circle-over} and
\eqref{cor-ho-inv-of-simpl-algebras}.

\subsection{Notation and conventions}\label{subsec-notations}

We write $X\backslash\cat{C}$ and $\cat{C}/X$ for the categories of
objects under and over a given object $X$ (the ``comma categories'').  
We write $\cat{D}^\cat{C}$ or $\Func(\cat{C},\cat{D})$ for the
category of functors from $\cat{C}$ to $\cat{D}$.

If $X$ and $Y$ are algebras over some monad $T$, we adopt the
convention of writing
$\displaystyle X\coprod^{\alg{T}}Y$ or $X\amalg^T Y$ for the coproduct
of $X$ and $Y$ in the category of $T$-algebras.  An undecorated
coproduct symbol means one taken in some underlying category,
which typically is sets or simplicial sets.

We write $\Set$ for the category of sets, and $\sSet$ for the category
of simplicial sets; it is often convienient to regard $\Set\subset
\sSet$ as the full subcategory of \emph{discrete} simplical sets.
Generally, we write $s\cat{C}$ for the category of 
simplicial objects in $\cat{C}$.  The diagonal functor $\diag\colon
s(s\cat{C})\ra s\cat{C}$ sends $\{Y_{p,q}\}\mapsto \{Y_{n,n}\}$.  We
often use the diagonal principle
\cite[IV.1.7]{goerss-jardine-simplicial-book}, which says that if
$f\colon X\ra Y$ is a 
morphism in $s(s\Set)$ (i.e., of bisimplicial sets) such that 
$f_{p,*}\colon X_{p,*}\ra Y_{p,*}$ is a weak equivalence of simplicial
sets for every $p\geq0$, then $\diag(f)$ is a weak equivalence of
simplicial sets.

\subsection{Acknowledgments}

The author would like to thank Paul Goerss for conversations which
improved the paper.  The author would also
like to thank Haynes Miller for suggesting an improved title.

\subsection{Organization of the paper}

In \S\ref{sec-proper-model-cats} we describe the notion of proper
model categories and prove some key properties; we also give several
examples of categories of simplicial algebras which are not proper.
In \S\S\ref{sec-functors-on-sets} and \ref{sec-theories-algebras-bimod}
we establish what we need for algebraic theories and their algebras
over sets and simplicial sets.  In the approach we take, algebraic
theories are simply monads over sets (or graded sets) which
commute with filtered colimits.  We also establish the notion of a
bimodule between theories, and identify them with a certain class of
functors between categories of algebras.  In
\S\S\ref{sec-functors-commuting-with-prod} and \ref{sec-s-free-maps} we
carry out some preparations needed for \S\ref{sec-ho-theory-of-alg},
in which we describe the model category structure on categories of
simplicial algebras, and for \S\ref{sec-ho-invariance-prop}, in which
we show that the homotopy theory of algebras over a theory is a weak
homotopy invariant of the theory, and that cofibrant right modules
over a theory preserve all weak equivalences.  In
\S\ref{sec-proper-criterion} we establish a criterion for a category
of simplicial algebras to be proper, by generalizing an argument of
Dwyer and Kan \cite{dwyer-kan-simplicial-localizations}.  In
\S\ref{sec-free-theories-and-trees} we give a description of free
theories using trees, which is then used in
\S\ref{sec-cofs-of-theories} to show that a cofibrant theory gives
rise to a proper model category of algebras.  We give proofs of
Theorems A--D in \S\ref{sec-proofs-of-theorems}.

\section{Proper model categories}\label{sec-proper-model-cats}

By \dfn{model category}, we mean a closed model category in the sense
of Quillen \cite{homotopical-algebra}, \cite{quillen-ratl-homotopy}.
(See also \cite{hovey-model-categories}, who defines model categories
with a slightly stronger set of axioms than Quillen.  However,
everything in this section holds under Quillen's axioms.)  
We write $\Ho \cat{M}$ for the category obtained by formally inverting
the weak equivalences in a model category $\cat{M}$.

\subsection{Definition of properness}\label{subsec-proper-def}

We recall the notion of a proper model category.

\begin{defn}
A  model category $\cat{M}$ is  \dfn{left proper} if
for each pushout square in $\cat{M}$ of the form
$$\xymatrix{
{A} \ar[r]^i \ar[d]_f 
& {B} \ar[d]^g
\\
{C} \ar[r]
& {D}
}$$
in which $i$ is a cofibration and $f$ is weak equivalence, the map
$g$ is weak equivalence.

Similarly, $\cat{M}$ is \dfn{right proper} if it satisfies the 
dual property involving pullback squares, fibrations, and weak
equivalences.

A  model category $\cat{M}$ is \dfn{proper} if it is both left
proper and right proper.
\end{defn}

\subsection{Under- and over-categories and properness}

Properness is most naturally understood as a statement about
``families'' of model
categories which are parameterized by the objects of a fixed model
category. 

We say that a pair of adjoint functors $L\colon
\cat{M}\rightleftarrows\cat{N}\noloc R$ between model categories is a
\dfn{Quillen pair} if the left 
adjoint $L$ takes cofibrations to cofibrations and the right adjoint
$R$ takes
fibrations to fibrations.  The pair forms a \dfn{Quillen equivalence}
if, in addition, for each cofibrant object $X$ in $\cat{M}$ and
fibrant object $Y$ 
in $\cat{N}$, a map $LX\ra Y\in\cat{N}$ is a weak equivalence if and
only if its adjoint $X\ra RY\in\cat{M}$ is.  
\begin{prop}
A Quillen pair as above gives rise to a derived adjoint pair
$\Ho\cat{M}\rightleftarrows \Ho\cat{N}$.  Furthermore, the derived
pair is an equivalence if and only if the Quillen pair is a Quillen
equivalence. 
\end{prop}
\begin{proof}
See \cite{abstract-homotopy} or \cite[1.3.10 and
1.3.13]{hovey-model-categories}. 
\end{proof}

Recall that given an object $X$ in a model category $\cat{M}$ the
categories $X\backslash\cat{M}$ and $\cat{M}/X$ of objects under and
over  $X$ are naturally equipped with model
category structures, in which the fibrations, cofibrations, and weak
equivalences are inherited from $\cat{M}$.  Furthermore, given a map
$f\colon X\ra Y$ in $\cat{M}$, the induced adjoint functor pairs
$$Y\amalg_X{-}\colon X\backslash\cat{M} \rightleftarrows
Y\backslash\cat{M}\noloc f^*
\qquad\text{and}\qquad f_*\colon \cat{M}/X \rightleftarrows \cat{M}/Y
\noloc X\times_Y{-}$$
are Quillen pairs.  We note that
\begin{prop}\label{prop-overcat-quillen-pair-equiv}
Let $\cat{M}$ be a model category, and suppose $f\colon X\ra Y\in
\cat{M}$.  Then the following are equivalent:
\begin{enumerate}
\item [(1)] The pair $X\backslash\cat{M}\rightleftarrows
Y\backslash\cat{M}$ (resp.\ $\cat{M}/X\rightleftarrows \cat{M}/Y$) is
a Quillen equivalence. 

\item [(2)]  The pushout (resp.\ pullback) of $f$ along any cofibration
(resp.\ fibration) in $\cat{M}$ is a weak equivalence.
\end{enumerate}
A necessary condition for (1) and (2) to hold is that $f$ be a weak
equivalence.  Sufficient conditions for (1) and (2) to hold are: that
$f$ be a trivial cofibration (resp.\ trivial fibration), \emph{or} that
$X$ and $Y$ be cofibrant (resp.\ fibrant) objects.
\end{prop}
\begin{proof}
We give the proof of the cofibration case, as the fibration case is
strictly dual.  Let $i\colon X\ra X'$ and $j\colon Y\ra Y'$.  Then a
map $g\colon X'\ra Y'\in X\backslash\cat{M}$ and its adjoint $g'\colon
X'\cup_X Y\ra Y'\in Y\backslash\cat{M}$ are related by $g=g'f'$, where
$f'$ is the pushout of $f$ along $i$.  If (1) holds and if $i$ is
a cofibration, we can construct 
$g'$ so that it is a weak equivalence to a fibrant object, and it then
follows from (1)
that $g$ and hence $f'$ are weak equivalences, giving (2).
Conversely, if (2) holds, then $g$ is a weak equivalence if and only
if $g'$ is, giving (1).

The necessary condition follows from considering the case when $i$ is
the identity map.  That $f$ being a trivial cofibration is sufficient
is clear; that $X$ and $Y$ being cofibrant is sufficient then follows
using \eqref{lemma-ken-brown} and the fact that Quillen equivalences
satisfy a 2 out of 3 property \cite[1.3.15]{hovey-model-categories}. 
\end{proof}

\begin{lemma}\label{lemma-ken-brown}
Let $f\colon X\ra Y\in\cat{M}$ be a weak equivalence between cofibrant
objects.  Then the exists a factorizaton $f=pi$ such that $p$ admits a
section $s$ and both $i$ and $s$ are trivial cofibrations.
\end{lemma}
\begin{proof}
See \cite[1.1.12]{hovey-model-categories}.
\end{proof}

Thus one has the following reformulation of the
notion of properness.
\begin{prop}\label{prop-undercat-char-of-proper}
A model category $\cat{M}$ is left (resp.\ right) proper if and only
if for every weak 
equivalence $f\colon X\ra Y$ in $\cat{M}$, the induced adjoint functor
pair $X\backslash\cat{M} \rightleftarrows Y\backslash\cat{M}$ (resp.\
$\cat{M}/X \rightleftarrows \cat{M}/Y$) is a 
Quillen equivalence.
\end{prop}

\begin{rem}\label{rem-properness-conditions}
\forcepar
\begin{enumerate}
\item [(i)] \label{rem-easy-properness} If all objects in $\cat{M}$
are cofibrant (resp.\ fibrant) then $\cat{M}$ is left (resp.\ right)
proper, because of the sufficiency condition of
\eqref{prop-overcat-quillen-pair-equiv}. 

\item [(ii)] Note that if $\cat{M}$ is left or right proper, then so
are all comma categories $X\backslash\cat{M}$ and $\cat{M}/X$.
\end{enumerate}
\end{rem}

We should note that proper model categories have good theories of
homotopy cartesian and cocartesian squares; this topic is treated in
detail in
\cite[II.8]{goerss-jardine-simplicial-book}.

\subsection{Examples of proper model categories}

The categories of simplicial sets and of topological spaces are
examples of proper model categories; for simplicial sets a proof is
given in \cite[II.8.6]{goerss-jardine-simplicial-book}.  
We will observe in \eqref{thm-closed-model-cat} that all categories of
simplicial algebras are \emph{right} proper.

In certain cases, model categories of simplical algebras (as defined
in \S\ref{sec-ho-theory-of-alg}) are known to
be left proper (and hence proper).  These examples include simplical
objects in: all
abelian categories, 
commutative monoids,  
monoids, simplicial categories with a fixed set of objects
\cite{dwyer-kan-simplicial-localizations}, 
commutative algebras over any commutative ring $R$, and associative
algebras over a field.
\margnote{
I'm not sure about simplicial groups.  Or simplicial lie algebras over
a field.}

\subsection{Examples of improper model categories}
\label{sec-improper-model-cat}

Not every category of simplicial algebras is \emph{left} proper.  We
offer two examples in which left properness fails.  In 
both cases, left properness is shown to fail by observing that it
fails in the simplest case: the functor which takes an simplicial
algebra $X$ to the 
coproduct of 
$X$ with a free algebra on one generator does not in general take weak
equivalences to weak equivalences (cf.\
\ref{thm-properness-iff-coproduct} and
\ref{rem-singleton-properness-criterion}).  
It
should be apparent that many other such examples could be constructed,
and that failure of left properness is a ``generic'' property of
categories of simplicial algebras. 
\begin{exam}
Let $T$ be the theory of associative algebras over a commutative ring
$R$.  If $R$ has $\Tor$-dimension greater than $0$, then simplicial 
$T$-algebras is not a left proper model category.  (This example
was pointed out to me by Paul Goerss.)

If $A$ is an associative $R$-algebra, the algebra $A\langle
x\rangle$ obtained by adjoining one free generator has the form
$$A\langle x\rangle = \bigoplus_{n\geq1} A^{\otimes n},$$
where the tensor product is taken over $R$.  (The $n$-fold tensor
product in this sum corresponds to all expressions in $A\langle
x\rangle$ of the form $a_1 xa_2 xa_3 x \dots x a_n$, with $a_i\in A$.)
Any $R$-algebra $A$ is weakly equivalent to a simplicial $R$-module
$B$ which is degreewise flat over $R$, by taking a free resolution.
Thus, if there exists an algebra $A$ such that $\Tor^R_i(A,A)\neq 0$
for some $i>0$ (e.g., $A=R\oplus M\oplus N$ with $\Tor_1^R(M,N)\neq0$
and with trivial product on $M\oplus N$), then $A\langle x\rangle$ 
is not weakly equivalent to $B\langle x\rangle$.  
\end{exam}
\begin{exam}
Let $C$ denote the theory of augmented commutative $R$-algebras.  Any such
algebra $A$ has an augmentation ideal $I(A)$.  Thus, this category is
equivalent to the category of non-unital commutative $R$-algebras, by the
functor sending $A\mapsto I(A)$.
Let $C_n$
for $n\geq 1$ denote the theory of augmented commutative $R$-algebras
with the additional property that $I(A)^n=0$.  The category of
simplicial algebras over $C_n$ is not proper for $n\geq 3$.

We give the proof in the case $n=3$; the general case is no more
difficult.  In this case, if $A$ is a 
$C_3$-algebra with augmentation ideal $I=I(A)$, and $A\langle
x\rangle$ is the $C_3$-algebra obtained by adjoining one 
free generator to $A$, we have
$$A\langle x\rangle \approx A[x]/(I,x)^3 \approx
A \oplus (A/I^2)x \oplus (A/I)x^2.$$
Since $A/I=R$ and $A/I^2 = R\oplus I/I^2$, the functor
$A\mapsto A\langle x\rangle$ is a direct sum of the 
identity functor, two copies of the functor with constant value $R$, and the
indecomposables functor $I\mapsto I/I^2$.  The indecomposables functor
on $C_3$-algebras has non-trivial 
higher derived functors, and hence if $B$ is a
free simplicial resolution of $A$, then $A\langle x\rangle$ will not
in general be weakly equivalent to $B\langle x\rangle$.  (A specific
example where this occurs is $A=R[y]/y^2$.)
\end{exam}

We also note for the record that there are examples of model
categories which are not right proper; the first was given by Quillen
\cite[II.2.9]{quillen-ratl-homotopy}.  Here is a more typical example.
Consider the
category of simplicial sets, equipped with a model category structure
in which weak equivalences are rational homology isomorphisms, and
cofibrations are inclusions; this is an example of Bousfield's
localization model category structure
\cite{bousfield-localization-spaces}.  Then one can form a pull-back
square of the form 
$$\xymatrix{
{K(\Q/\Z,0)} \ar[r]_-g \ar[d]
& {C} \ar[d]^p
\\
{K(\Z,1)} \ar[r]^f 
& {K(\Q,1)} 
}$$
in which $p$ is a rational fibration from a contractible space $C$ and $f$
is a rational homology 
isomorphism, but $g$ is not a rational homology isomorphism.

\section{Functors on sets}\label{sec-functors-on-sets}

In this section we characterize functors between
categories of sets (and more generally, graded sets) which commute
with filtered colimits.  This is a prerequisite to our approach to
theories in \S\ref{sec-theories-algebras-bimod}.  

\subsection{Reflexive coequalizers}

A \dfn{reflexive pair} in a category is a diagram consisting of a pair
of maps $f,g\colon X\ra 
Y$ together with a map $s\colon Y\ra X$ (called a \dfn{reflection})
such that $fs = 1_Y = gs$.  The colimit
of such a diagram is the same as the coequalizer of the pair
$f,g$; we call it a \dfn{reflexive coequalizer}.  We record the
following elementary but useful fact.
\begin{prop}
In $\Set$ (the category of sets), reflexive coequalizers commute with
finite products. 
\end{prop}

\subsection{Functors from finite sets}

Let $\fsets\subset\Set$ denote a fixed skeleton of the
full subcategory of finite sets.  It will be convenient to identify 
$\ob\fsets$ with $\N$, and to write $\fset{n}=\{1,\dots,n\}\in
\ob\fsets$ for the distinguished copy of the $n$-element set.  We
write $X^{\fset{n}}=\hom_\Set(\fset{n},X)$.  

Let $r\colon \Func(\Set,\Set) \ra \Set^\fsets$ denote the restriction
functor; it takes an endofunctor on sets to its restriction $\fsets\ra
\Set$.  This functor admits a left adjoint $\iota\colon \Set^\fsets\ra
\Func(\Set,\Set)$, which associates to each $A\in \Set^\fsets$ its
\emph{left Kan extension} $\iota A\colon \Set\ra \Set$ along the full embedding
$\fsets\subset\Set$.  This can be presented as a
reflexive coequalizer:
\begin{equation}\label{eq-defining-coequalizer-for-action}
\coprod_{\fset{p}\ra\fset{q}} A(\fset{p})\times X^{\fset{q}}
\rightrightarrows \coprod_{\fset{n}} A(\fset{n})\times X^{\fset{n}}
\ra (\iota A)(X),
\end{equation}
where $X\in \Set$ and $A\in \Set^\fsets$.
Because $\fsets\subset\Set$ is a \emph{full} subcategory, one sees
that $A\ra r\iota A$ is an isomorphism for all $A$ in $\Set^\fsets$,
(that is, $(\iota 
A)(\fset{n}) \approx A(\fset{n})$) and hence that $\iota$
identifies $\Set^\fsets$ up to equivalence as a full subcategory of
the category of all functors.

Let $\filtofunc(\cat{C},\cat{D})\supset \filtfunc(\cat{C},\cat{D})$ denote the
full subcategories of $\Func(\cat{C},\cat{D})$ consisting of those functors
which commute respectively: with filtered colimits;
with filtered 
colimits and reflexive coequalizers.  Note that both
subcategories are
closed under composition of functors.
\begin{prop}\label{prop-char-of-filt-functors-set}
There is a factorization $\iota\colon \Set^\fsets\ra \Func(\Set,\Set)$
into
$$\Set^\fsets \xra{\iota_1} \filtfunc(\Set,\Set)
\xra{\iota_2}\filtofunc(\Set,\Set) \subset \Func(\Set,\Set),$$
and $\iota_1$ and $\iota_2$ are equivalences of categories.
\end{prop}
\begin{proof}
Let $A\in \Set^\fsets$.  The formula
\eqref{eq-defining-coequalizer-for-action} for $\iota A$ given above shows that
$\iota A(X)$ is computed from $A$ and $X$ using only colimits and finite
products, and both of these commute with filtered colimits and
reflexive coequalizers.  Hence $\iota$ factors through a functor $\iota_1$.  

To show that $\iota_1$ and $\iota_2$ are equivalences, it suffices to
show that
if $F\in \filtofunc(\Set,\Set)$, then $\eta_F\colon \iota rF\ra F$ is an
isomorphism.  In fact, $\eta_F$ is clearly an isomorphism when evaluated
at any finite set, and the result follows from the fact that every set
is a filtered colimit of its finite subsets.
\end{proof}

As an immediate consequence of \eqref{prop-char-of-filt-functors-set}
we see that $\Set^\fsets$ admits the structure of a monoidal category,
which corresponds via $\iota$ to composition of functors.  We will
denote this monoidal structure by $A\circ B$, for $A,B\in
\Set^\fsets$, so that $\iota(A\circ B)\approx \iota A\circ \iota B$.
The unit corresponds to $I$, defined by 
$I(n)=\hom_{\fsets}(\fset{1},\fset{n})=\fset{n}$.  Given $A,B\in
\Set^\fsets$ and $\fset{m}\in \fsets$, a formula for $(A\circ
B)(\fset{m})$ can be derived by inserting $B(\fset{m})$ for $X$ in
\eqref{eq-defining-coequalizer-for-action}.

\subsection{Graded sets}

Let $\indset{I}$ be a set.  An \dfn{$\indset{I}$-graded set} is a collection
$(X_i)_{i\in \indset{I}}$ of sets, and a morphism of such is
collection of maps respecting the grading.  The category of
$\indset{I}$-graded sets is denoted $\Set^\indset{I}$.  An
$\indset{I}$-graded set is said to be \dfn{finite} if $\coprod_i X_i$
is a finite set.

We write
$\fsets/\indset{I}$ for the category whose objects are functions $f\colon
\fset{n}\ra \indset{I}$, $\fset{n}\in\fsets$, and whose morphisms are
commuting triangles.  A function $f\colon X\ra \indset{I}$ of sets
naturally gives rise to 
an $\indset{I}$-graded set $(f^{-1}(i))_{i\in\indset{I}}$, giving an
inclusion 
functor  
$\fsets/\indset{I}\ra \Set^\indset{I}$ which is equivalent to the
inclusion of the full subcategory of \emph{finite} $\indset{I}$-graded
sets.  

Given sets $\indset{I}$ and $\indset{J}$, let $\fsetsdind{I}{J} =
\indset{J}\times \fsets/\indset{I}$.  Objects in this category
are pairs $(j\in\indset{J}, f\colon\fset{n}\ra\indset{I})$; a morphism
$(j,f)\ra (j',f')$ is defined only if $j=j'$, in which case it
consists of a map $f\ra f'\in \fsets/\indset{I}$.  We write
$\N(\indset{I},\indset{J})=\ob\fsetsdind{I}{J}$.  
Then
$\Set^\fsetsdind{I}{J} = (\Set^\indset{J})^{\fsets/\indset{I}}$ is
equivalent to the category of functors 
from \emph{finite} $\indset{I}$-graded sets to $\indset{J}$-graded
sets, giving rise to a restriction functor $r\colon
\Func(\Set^\indset{I},\Set^\indset{J})\ra \Set^\fsetsdind{I}{J}$.
This functor admits a left adjoint $\iota \colon \Set^\fsetsdind{I}{J}\ra
\Func(\Set^\indset{I},\Set^\indset{J})$, which associates to each
$A\in \Set^\fsetsdind{I}{J}=(\Set^\indset{J})^{\fsets/\indset{I}}$ its
\emph{left Kan extension} $\iota A\colon 
\Set^\indset{I}\ra \Set^\indset{J}$ along the full embedding
$\fsets/\indset{I}\subset \Set^\indset{I}$.  There is a reflexive
coequalizer formula:
\begin{equation}\label{eq-defining-coequalizer-for-gr-action}
\coprod_{p\ra q\in \fsets/\indset{I}} A(j,p)\times X^q
\rightrightarrows \coprod_{f\in\ob\fsets/\indset{I}} A(j,f)\times X^f \ra
(\iota A)(X)_j,\qquad j\in\indset{J},
\end{equation}
where for $f\colon \fset{n}\ra \indset{I}\in
\ob\fsets/\indset{I}\subset \Set^\indset{I}$ and 
$X\in \Set^\indset{I}$ 
we write $X^f = \hom_{\Set^\indset{I}}(f,X) = \prod_{k\in \fset{n}}
X_{f(k)}\in\Set$. 

Since $\fsets/\indset{I}\ra\Set^\indset{I}$ is full, we have that
$A\approx r\iota A$, so that $(\iota A)(K)_j \approx A(j,K)$ for
$j\in\indset{J}$ and $K\in\fsets/\indset{I}\subset\Set^\indset{I}$.
We take advantage of 
this fact to write $A(K)\in \Set^\indset{J}$ for the
$\indset{J}$-graded set $A(\mbox{--},K)$.

Considerations identical to those which gave
\eqref{prop-char-of-filt-functors-set} give
\begin{prop}\label{prop-char-of-filt-functors-gr-set}
There is a factorization $\iota\colon \Set^\fsetsdind{I}{J}\ra
\Func(\Set^\indset{I},\Set^\indset{J})$ 
into
$$\Set^\fsetsdind{I}{J} \xra{\iota_1}
\filtfunc(\Set^\indset{I},\Set^\indset{J}) 
\xra{\iota_2}\filtofunc(\Set^\indset{I},\Set^\indset{J}) \subset
\Func(\Set^\indset{I},\Set^\indset{J}),$$ 
and $\iota_1$ and $\iota_2$ are equivalences of categories.
\end{prop}

An immediate consequence of \eqref{prop-char-of-filt-functors-gr-set}
is the existence of pairings $\mbox{--}\circ\mbox{--}\colon
\Set^\fsetsdind{J}{K}\times \Set^\fsetsdind{I}{J}\ra
\Set^\fsetsdind{I}{K}$ which correspond via $\iota$ to composition of
functors.  
Let $\fsetsind{I}=\fsetsdind{I}{I}$ and $\N(\indset{I}) =
\ob\fsetsind{I}$.  Then  
$\Set^\fsetsind{I}$ is a monoidal category, and in fact is a full monoidal
subcategory of the category of endofunctors of $\Set^\indset{I}$.  
The unit object $I$ is defined by $I(i,f\colon \fset{n}\ra \indset{I})
= 
f^{-1}(i)$.  

\subsection{Free series}\label{subsec-free-series}

Let $\indset{I}$ and $\indset{J}$ be sets.  Recall that
$\N(\indset{I},\indset{J}) = \ob \fsetsdind{I}{J}$.  The forgetful
functor $\Set^\fsetsdind{I}{J}\ra 
\Set^{\N(\indset{I},\indset{J})}$ admits a left adjoint $\freeser
\colon \Set^{\N(\indset{I},\indset{J})}\ra \Set^\fsetsdind{I}{J}$
called the \dfn{free series} functor.  One easily checks the formula
$$\freeser A\approx \coprod_{K\in\N(\indset{I},\indset{J})} A(K)\times
I^K,$$
where $I^K\in \Set^\fsetsdind{I}{J}$ is defined by $L\mapsto
I^K(L)=\hom_{\fsetsdind{I}{J}} (K,L)$.  For $X\in \Set^\indset{I}$ we have 
$$(\iota\freeser A)(X)_j \approx
\coprod_{f\colon\fset{n}\ra\indset{I}\in\fsets/\indset{I}}
A(f)_j\times X^f,$$
where $X^f\in\Set$ is as in
\eqref{eq-defining-coequalizer-for-gr-action}.  In the case when
$\indset{I}=\indset{J}$ is a singleton, these formulas reduce to
$\freeser A \approx \coprod_n A(n)\times I^n$ and
$(\iota \freeser A)(X) = \coprod_n A(n)\times X^n$.

\subsection{Simplicial objects}

By prolongation, we obtain a functor $\sSet^\fsetsdind{I}{J}\ra
\Func(\sSet^\indset{I}, \sSet^\indset{J})$.  The image of this functor
is the full subcategory of \emph{simplicial objects} in the category of
functors which commute with filtered colimits and reflexive coequalizers.  
Formulas \eqref{eq-defining-coequalizer-for-action} and
\eqref{eq-defining-coequalizer-for-gr-action} still apply in this
case, where the objects are now \emph{graded} simplicial sets.

\section{Theories, algebras, and
bimodules}\label{sec-theories-algebras-bimod}

In this section, we define algebraic theories and their associated
algebra categories.  In our approach, we also consider multi-sorted
theories.  We also give some attention to \emph{bimodules} of
theories, which give rise to a large class of functors between
categories of algebras, and will play an important role in
\S\S\ref{sec-ho-invariance-prop} and \ref{sec-cofs-of-theories}.
The definitions of theories and algebras that we give appear quite
different than the notions of algebraic theories and their models as
in \cite{lawvere-functorial-semantics}, where a theory is defined to
be a category with finite products (see the nice treatment in
\cite[Ch.\ 3]{borceux-handbook-cat-alg-2} for this).  However, our
categories 
of ``algebras'' are the same as the categories of ``models'', as
we note 
below \eqref{subsec-algebras}.  Our formulation is one of those used
by Boardman and Vogt in a topological context
\cite{boardman-vogt-homotopy-invariant-structures} (they write
``theories with colours'' for what we call ``multi-sorted theories'').  
It is also close to that given by 
Schwede \cite{schwede-stable-ho-alg-theor}, although the theories he
considers are \emph{pointed}.  Lawvere's original formulation is also
used in \cite{boardman-vogt-homotopy-invariant-structures}, and
is used in a crucial way by Badzioch \cite{badzioch-alg-theories}.

In what follows, we make repeated use of the identifications
of
$\Set^\fsetsdind{I}{J}$ and $\sSet^\fsetsdind{I}{J}$ as full
subcategories of the respective functor categories, and we omit use of
the $\iota$ symbol of \S\ref{sec-functors-on-sets}; thus we write
$A(X)$ where before we had $(\iota A)(X)$.

\subsection{Theories}\label{subsec-theories}

Let $\indset{I}$ be a set, and recall that $\Set^\fsetsind{I}$ is a
monoidal category, equivalent to a full monoidal subcategory of
$\Func(\Set^\indset{I}, \Set^\indset{I})$.  We define an
\dfn{$\indset{I}$-sorted theory}, or more simply a
\dfn{theory}, to be a monoid object $T$ in $\Set^\fsetsind{I}$.  That
is, $T\in\Set^\fsetsind{I}$ is equipped with maps $\mu_T\colon T\circ T\ra T$ and
$\eta_T\colon I\ra T$ satisfying the usual axioms for a monoid.  
From
\eqref{prop-char-of-filt-functors-gr-set}, we see that
$\indset{I}$-sorted theories are
essentially the same as monads on $\Set^\indset{I}$ which commute with
filtered colimits.

We write $\theories(\indset{I})$ for the category of
$\indset{I}$-sorted theories over sets.  

\subsection{Algebras over a theory}\label{subsec-algebras}

An \dfn{algebra} $X$ over an $\indset{I}$-sorted theory $T$ is an
algebra over the monad induced by $T$; that is, an algebra is an
object $X\in \Set^\indset{I}$ equipped with a map $\psi\colon T(X)\ra
X$ satisfying the usual axioms.  The category of $T$-algebras is
denoted $\alg{T}$.

Given a graded set $X$, the object $T(X)$ is
naturally a $T$-algebra, namely the \dfn{free $T$-algebra on $X$}.  

\begin{prop}\label{prop-algebras-complete-cocomplete}
Let $T$ be an $\indset{I}$-sorted theory.  The
category $\alg{T}$ is complete and cocomplete.  Limits, filtered
colimits and reflexive coequalizers are created in the underlying
category $\Set^\indset{I}$.  There exists an adjoint functor pair
$T\colon \Set^\indset{I}\rightleftarrows \alg{T}\noloc u$, where $u$
is the forgetful functor, and $T$ is called the free $T$-algebra
functor. 
\end{prop}
\begin{proof}
That limits, filtered colimits, and reflexive coequalizers exist and
are created in $\Set^\indset{I}$ is immediate from
\eqref{prop-char-of-filt-functors-gr-set}.  That the free algebra
functor is left adjoint is a standard property
of monads \cite[4.1.4]{borceux-handbook-cat-alg-2}.  Existence of
colimits follows from \cite[4.3.6]{borceux-handbook-cat-alg-2}; or
note that  
colimits of a diagram $\alpha\mapsto X_\alpha \colon \cat{A}\ra
\alg{T}$ can be constructed 
explicitly as the 
reflexive coequalizer in $\alg{T}$ of
$T(\colim^{\Set^\indset{I}}_{\cat{A}} TX_\alpha) 
\rightrightarrows T(\colim^{\Set^\indset{I}}_{\cat{A}} X_\alpha)$, the
top map being induced by the inclusions $X_\alpha
\ra \colim^{\Set^\indset{I}}_{\cat{A}} TX_\alpha$ and the top map
being induced by the algebra structure maps $TX_\alpha\ra
X_\alpha$.
\end{proof}

According to our definition, a single-sorted theory $T$ corresponds,
via \eqref{prop-char-of-filt-functors-set}, precisely to a monad on
sets which commutes with filtered colimits.  By
\cite[4.6.2]{borceux-handbook-cat-alg-2}, categories of algebras over
such monads (which we have called ``algebras over a theory'') are
exactly those which are equivalent to categories of ``models of an
algebraic theory'' in the classical sense.  See also
\cite[Prop. 2.30]{boardman-vogt-homotopy-invariant-structures}.  

Finally, note that the category $\theories(\indset{I})$ of
$\indset{I}$-sorted 
theories is itself an example of a category of algebras over a certain
$\N(\indset{I})$-sorted theory, namely the theory of
$\indset{I}$-sorted theories.  This is because the forgetful functor
$\theories(\indset{I}) \ra \Set^{\N(\indset{I})}$ admits a left
adjoint, and is\margnote{need to say more here} monadic.
Thus
\eqref{prop-algebras-complete-cocomplete} shows that the category of
such theories is complete and cocomplete, and that there exist free
theories.  We will consider an explicit construction of free theories
in \S\ref{sec-free-theories-and-trees}.

\subsection{Bimodules}\label{subsec-bimodules-def}

Given $S\in \theories(\indset{I})$ and $T\in \theories(\indset{J})$, a
\dfn{$T,S$-bimodule} is an object $M\in \Set^\fsetsdind{I}{J}$
equipped with actions $T\circ M\ra M$ and $M\circ S\ra M$, which are
associative and unital and which 
commute with each other.  Let $\bimod{T,S}$ denote the category of
bimodules.  A \dfn{right $S$-module} is an $I,S$-bimodule
and a \dfn{left $T$-module} is a $T,I$-bimodule.

Given an $S$-algebra $X$, let $M\circ_S X$
denote the coequalizer of the following reflexive pair in $\alg{T}$
(which can be 
computed in graded sets by \eqref{prop-algebras-complete-cocomplete}):
$$M(S(X)) \rightrightarrows M(X) \ra M\circ_S X.$$
This gives rise to a functor $\iota\colon
\bimod{T,S}\ra\Func(\alg{S},\alg{T})$.   (Warning: this is not the
$\iota$ used in \S\ref{sec-functors-on-sets}.) 
Note that if $K\in\fsets/\indset{I}\subset \Set^\indset{I}$, then
$M\circ_S S(K)\approx M(K)$ as objects of $\Set^\indset{J}$; that is,
$M(K)$ is the value of $M\circ_S\mbox{--}$ on 
the free $S$-algebra generated by $K$.
\begin{prop}\label{prop-char-of-filt-functors-algebras}
Let $S$ and $T$ be $\indset{I}$- and $\indset{J}$-sorted theories over
sets.  The functor $\iota\colon \bimod{T,S}\ra \Func(\alg{S},\alg{T})$ defined
above factors through,
and induces an equivalence with, the full subcategory
$\filtfunc(\alg{S},\alg{T})$ of functors which commute with filtered
colimits and reflexive coequalizers.
\end{prop}
\begin{proof}
It is clear using \eqref{prop-char-of-filt-functors-gr-set} that
$\iota$ factors through the subcategory.  It remains to show that
$\iota$ is an equivalence.

Let $\alg{S}^\fgf\subset \alg{S}$ denote the full subcategory of
\emph{finitely generated free} $S$-algebras; every object in this
subcategory is isomorphic to $S(K)$ for some $K\in\fsets/\indset{I}$.
Consider the sequence of functors
$$\bimod{T,S} \xra{\iota} \filtfunc(\alg{S},\alg{T}) \subset \Func(\alg{S},
\alg{T}) \ra \Func(\alg{S}^\fgf, \alg{T});$$
the right-hand arrow is the one induced by restriction of functors to
the subcategory.  The result will follow when we show that the
composites $\alpha\colon \filtfunc(\alg{S},\alg{T}) \ra
\Func(\alg{S}^\fgf, \alg{T})$ and $\beta\colon \bimod{T,S}\ra
\Func(\alg{S}^\fgf, \alg{T})$ are equivalences.

To see that $\alpha$ is an equivalence, observe that every $S$-algebra
is a coequalizer of a reflexive diagram of free algebras, and that
every free $S$-algebra is a filtered colimit of finitely generated
free algebras.  
Thus every functor $\alg{S}\ra \alg{T}$ which commutes
with filtered colimits and reflexive coequalizers is
determined up to unique isomorphism by its restriction to the subcategory
of finitely generated free 
algebras, and natural transformation between such functors are
uniquely determined by this restriction.  Any
functor $\alg{S}^\fgf\ra \alg{T}$ extends to an element of
$\filtfunc(\alg{S},\alg{T})$ by a left Kan extension construction, and
therefore this construction gives the inverse to $\alpha$.

We now show that $\beta$ is an equivalence.  Explicitly, $\beta$ sends
$M\in \bimod{T,S}$ to the functor $G\colon X\mapsto M\circ_S X$; note
that if $X\approx S(K)$, then $G(X)\approx M(K)$.  We will construct
an inverse
$\gamma\colon \Func(\alg{S}^\fgf, \alg{T})\ra \bimod{T,S}$.
Given $F\colon \alg{S}^\fgf \ra \alg{T}$, define $N\in
\Set^\fsetsdind{I}{J}$ by $N(K) = F(S(K))$; recall that under the equivalence
$\Set^\fsetsdind{I}{J} \approx \filtfunc(\Set^\indset{I},
\Set^\indset{J})$, the object $N$ corresponds to a functor $\Set^\indset{I}\ra
\Set^\indset{J}$ also denoted by $N$, and there is a map
$N(X)\ra F(S(X))$ natural in $X\in \Set^\indset{I}$.  
Give $N$ the structure of a left $T$-module by 
$$(T\circ N)(K)=
T(F(S(K))) \ra F(S(K))= N(K),$$ 
using the fact that $F$ takes values in
$T$-algebras.  Give $N$ the structure of a right $S$-module by
$$(N\circ S)(K)\approx N(S(K)) \ra F(S(S(K))) \xra{F(\mu_K)} F(S(K)),$$
using the $S$-algebra structure of $S(K)$.  It follows that $N$ is a
$T,S$-bimodule, that our construction $F\mapsto \gamma F = N$ is a
functor from functors to bimodules, and that $\beta\gamma\approx \id$
and $\gamma\beta\approx \id$ as desired.

\end{proof}

\begin{rem}\label{rem-circ-on-right-has-bimod}
Let $S\in\theories(\indset{I})$ and $T\in \theories(\indset{J})$.
Then the category $\bimod{T,S}$ is a category of algebras over a
certain $\N(\indset{I},\indset{J})$-sorted theory
$B_{T,S}$; this is because bimodules are simply algebras over the
monad $A\mapsto T\circ A\circ S$ on $\N(\indset{I},\indset{J})$-graded
sets, which commutes with filtered colimits.

Suppose that $X\in\alg{S}$, and consider the functor
$\bimod{T,S}\ra\alg{T}$ given by $M\mapsto M\circ_SX$.  This functor
commutes with filtered colimits and reflexive coequalizers by
construction, and thus by \eqref{prop-char-of-filt-functors-algebras}
it is represented by a certain $T,B_{T,S}$-bimodule $N_X$ (whose
underlying set is graded by
$\N(\N(\indset{I},\indset{J}),\indset{J})$!)  We will 
need this observation in \S\ref{sec-ho-invariance-prop}. 
\end{rem}

Given a morphism $\phi\colon S\ra T$ of $\indset{I}$-sorted theories,
there is an evident restriction functor $\phi^*\colon \alg{T}\ra
\alg{S}$, which is the identity on underlying graded sets.
\begin{prop}\label{prop-adjoint-pair-between-algebras}
The restriction functor $\phi^*$ admits a left adjoint functor
$\phi_*\colon 
\alg{S}\ra \alg{T}$, and $\phi_*X\approx T\circ_S X$.  
\end{prop}
\begin{proof}
Let $Y\in \alg{T}$.  Then $\hom_{\alg{T}}(\phi_*X,Y)$ is the equalizer
of $\hom_{\alg{T}}(TX,Y)\rightrightarrows \hom_{\alg{T}}(TSX,Y)$, or
equivalently of $\hom_{\Set^\indset{J}}(X,Y)\rightrightarrows
\hom_{\Set^\indset{J}}(SX,Y)$, where the two arrows send $f\colon X\ra
Y$ to $f(\psi_X)$ and $(\psi_Y)(\phi X)(Sf)$ respectively, where
$\psi_X\colon SX\ra X$ and $\psi_Y\colon TY\ra Y$ denote the algebra
structure maps. 
\end{proof}

\subsection{Undercategories and coproducts of
theories}\label{subsec-undercats} 

Let $T$ be an $\indset{I}$-sorted theory, and $X$ a $T$-algebra.  Define $T_X\in
\Set^\fsetsind{I}$ by $T_X(K) = T(K)\amalg^T X$.
\begin{prop}\label{prop-undercat-theory}
The object $T_X$ admits the structure of a theory, and there is an
equivalence of categories $\alg{T_X}\approx X\backslash\alg{T}$. 
\end{prop}
\begin{proof}
Define\margnote{this needs work} $I\ra T_X$ to be the evident map
$I(K)\approx K\ra 
T(K)\amalg^T X\approx T_X(K)$, and $T_X\circ T_X\ra T_X$
to be the evident 
map $(T_X\circ T_X)(K) \approx T_X(T_X(K)) \approx
T(T(K)\amalg^T X)\amalg^T X \ra T(K)\amalg^T X\approx
T_X(K)$.  Then $T_X$ is easily seen to be a theory, and the
evident functor $\alg{T_X}\ra X\backslash\alg{T}$ an equivalence.
\end{proof}

Given $X\in \Set^\indset{I}$ and $g\colon X\ra Y\in \Set^\indset{I}$,
there exist \dfn{endomorphism theories} $\Endth{X}$ and $\Endth{g}$,
with the property that $\hom_{\theories(\indset{I})}(T,\Endth{X})$ is
in bijective correspondence with the set of $T$-algebra structures on
$X$, and $\hom_{\theories(\indset{I})}(T,\Endth{g})$ is in bijective
correspondence with the set of pairs of $T$-algebra structures on $X$
and $Y$ which make $f$ a map of $T$-algebras.  They are given by the
formulas
$$\Endth{X}(j,f) = \hom_{\fsetsdind{I}{J}}(X^f,X_j),
$$
$$\Endth{g}(j,f) =
\hom_{\fsetsdind{I}{J}}(X^f,X_j)\times_{\hom_{\fsetsdind{I}{J}}(X^f,
Y_j)} \hom_{\fsetsdind{I}{J}}(Y^f,Y_j).$$  
Here $X^f$ is as in \eqref{eq-defining-coequalizer-for-gr-action}. 

Let $S$ and $T$ be two $\indset{I}$-indexed theories.  Then
$S\amalg^{\theories(\indset{I})} T$ denotes the coproduct in
$\theories(\indset{I})$.  
\begin{prop}\label{prop-alg-over-coprod-of-theories}
The category $\alg{(S\amalg^{\theories(\indset{I})} T)}$ has as
objects $X\in \Set^\indset{I}$ 
equipped with both an $S$-algebra and a $T$-algebra structure, and as
morphisms those maps which commute with both algebra structures.
\end{prop}
\begin{proof}
This is immediate from the existence of the endomorphism theories.  
\end{proof}

\subsection{Simplicial objects}\label{subsec-simplicial-alg}

We can similarly consider
\dfn{simplicial theories}, namely monoid objects in
$\sSet^\fsetsind{I}$; these are the same as simplicial objects in
$\theories(\indset{I})$, and we write $s\theories(\indset{I})$ for the
category of $\indset{I}$-simplicial 
theories.  

If $T$ is a simplicial theory, then by a $T$-algebra $X$ we mean a
\dfn{simplicial algebra}, namely an object of
$\sSet^\indset{I}$ which is an algebra over the monad induced by $T$.
Effectively, if $T=\{T_n\}$ is a simplicial theory, $X$ amounts to a
collection $\{X_n\}$ of objects in $\Set^\indset{I}$, such that each
$X_n$ is equipped with the structure of a $T_n$ algebra, together
with, for each $\delta\colon [m]\ra [n]\in\simpind$, a map
$X_\delta\colon X_n\ra
(T_\delta)^*X_m$ of $T_n$-algebras, with the conditions
$X_{\delta'}X_{\delta}=X_{\delta\delta'}$.  

Similarly, we have \dfn{simplicial bimodules}; these are objects $M$ in
$\sSet^\fsetsdind{I}{J}$ which are $T,S$-bimodules, or equivalently a
collection $M=\{M_n\}$ of objects in $\Set^\fsetsdind{I}{J}$ such that
each $M_n$ is a $T_n,S_n$-bimodule, together with simplicial operators
acting as above.

Propositions \eqref{prop-algebras-complete-cocomplete},
\eqref{prop-adjoint-pair-between-algebras},
\eqref{prop-undercat-theory}, and
\eqref{prop-alg-over-coprod-of-theories} carry over to the simplicial
setting: change $\Set$ to $\sSet$.  There is also a
simplicial analogue of \eqref{prop-char-of-filt-functors-algebras}.
We say a functor $F\colon\alg{S}\ra\alg{T}$ between categories of
simplicial algebras is \dfn{degreewise} if there exist functors
$F_n\colon \alg{S_n}\ra\alg{T_n}$ such that $F(X)_n\approx F_n(X_n)$
for $n\geq0$, together with natural transformations
$F_\delta(S_\delta)^*\colon 
F_n\ra (T_\delta)^*F_m$ for each $\delta\colon [m]\ra [n]\in \simpind$
satisfying the appropriate identities. 
\begin{prop}\label{prop-char-of-filt-functors-s-algebras}
Let $S$ and $T$ be $\indset{I}$ and $\indset{J}$-sorted simplicial
theories.  Then the functor $\iota\colon \bimod{T,S}\ra
\Func(\alg{S},\alg{T})$ factors through,
and induces an equivalence with, the full subcategory
$\dfiltfunc(\alg{S},\alg{T})$ of degreewise functors which in each
degree commute with filtered colimits and reflexive coequalizers.
\end{prop}
\begin{proof}
Apply \eqref{prop-char-of-filt-functors-algebras} in each simplicial
degree. 
\end{proof}

\section{Functors commuting with
products}\label{sec-functors-commuting-with-prod} 

Henceforward, we consider only \emph{simplicial} $\indset{I}$-sorted
theories and simplicial algebras over such, unless otherwise indicated.

Let $E\colon \sSet\ra \sSet$ be a functor.  Such a functor induces
a functor $\sSet^\indset{I}\ra \sSet^\indset{I}$, which we also
denote by $E$.  We say that $E$ \dfn{commutes with products} if
$E(1)\approx 1$ and for all $X,Y\in\sSet$ the natural map $E(X\times
Y)\ra EX\times EY$ is an isomorphism.  Note that if $f\colon
\fset{p}\ra \fset{n}\in\fsets\subset\Set$ and if we write $X^f\colon
X^{\fset{n}}\ra 
X^{\fset{p}}$ for the induced map on products, then $E(X^f)\approx
(EX)^f$.  

Suppose that in addition there is a natural transformation $\eta\colon
\Id\ra E$.  Then there exist natural maps $X\times EY\ra E(X\times
Y)$ for all $X,Y\in \sSet$; furthermore, these maps are coherent, in
the sense that both ways to obtain a map $X\times Y\times EZ\ra
E(X\times 
Y\times Z)$ are the same.  In particular, 
$E$ is a \emph{simplicial} functor.  

Let $F=\freeser A\in \sSet^\fsetsdind{I}{J}$ be a free series
\eqref{subsec-free-series} on some 
$A\in \sSet^{\N(\indset{I},\indset{J})}$, and let $X\in
\sSet^\indset{I}$.  The above discussion shows 
that there is an evident map $\alpha\colon F(E(X))\ra E(F(X))$ in
$\sSet^\indset{J}$; for instance, in the case when $\indset{I}=*$, the
map is 
$$\alpha\colon  F(E(X))\approx \coprod_{\fset{n}} A(\fset{n})\times
(EX)^{\fset{n}} 
\ra \coprod_{\fset{n}} E\left( A(\fset{n}) \times X^\fset{n} \right)
\ra 
E\left( \coprod_{\fset{n}} A(\fset{n})\times X^\fset{n} \right)
\approx E(F(X)),$$
induced by $(EX)^\fset{n}\approx E(X^\fset{n})$ and $X\times EY\ra
E(X\times Y)$.  

\begin{prop}\label{prop-prod-pres-map}
Let $E\colon \sSet\ra \sSet$ be a functor commuting with products
and equipped with a 
natural transformation $\eta\colon \Id\ra E$.  Then for all $X\in
\sSet^\indset{I}$ and $F\in \sSet^{\fsetsdind{I}{J}}$, 
there exist maps
$$\widetilde{\alpha}\colon F(E(X)) \ra E(F(X))$$
which are natural in $X$ and $F$, and which in the case that
$F=\freeser A$ is the map $\alpha$ described above.
Furthermore, the $\widetilde{\alpha}$ are the unique collection of
maps with this 
property. 
\end{prop}
\begin{proof}
For convenience, we write the proof only in the case $\indset{I}$ and
$\indset{J}$ are
singleton sets; the general case is only notationally more difficult. 

We first show that $\alpha$ is in fact a natural transformation
between functors defined on the full subcategory of free objects in
$\sSet^\fsets$. 
Consider a map $\freeser A\ra \freeser B$ between free objects.
This amounts to a
collection of maps $A(\fset{n})\ra \coprod_p B(\fset{p})\times
\fset{n}^\fset{p}$, $n\geq0$, and the induced map $\freeser A(X)\ra \freeser
B(X)$ factors
$$\coprod_n A(\fset{n})\times X^\fset{n} \ra \coprod_{n,p}
B(\fset{p})\times \fset{n}^\fset{p}\times X^\fset{n}\ra \coprod_p
B(\fset{p})\times X^\fset{p}.$$
To show that $\alpha$ commutes with these maps reduces to showing that
it commutes with $\fset{n}^\fset{p}\times X^\fset{n}\ra X^\fset{p}$,
which is clear.

Define $\widetilde{\alpha}=\alpha$ on the full subcategory of free
objects in $\sSet^\fsets$.  Since every object of $\sSet^\fsets$ is a 
reflexive coequalizer of a pair of free objects, the $\widetilde{\alpha}$
extend in a 
unique way to arbitrary objects in $\sSet^\fsets$.
\end{proof}

By the uniqueness property of \eqref{prop-prod-pres-map}, we see that 
$\widetilde{\alpha} \colon\Id(E(X))\ra E(\Id(X))$ is the identity, and
the two ways 
of getting 
maps $G(F(E(X)))\ra E(G(F(X)))$ coincide: $\widetilde{\alpha}_{GF}=(
\widetilde{\alpha}_G 
F)(G\widetilde{\alpha}_F)$.

\begin{cor}\label{cor-prod-pres-maps-on-alg}
Let $T$ be a (possibly multi-sorted) simplicial theory.
A product preserving functor $E\colon \sSet\ra \sSet$ equipped with a
natural transformation $\eta\colon \Id\ra E$ lifts in a
natural way to a functor $E\colon \alg{T}\ra \alg{T}$.
Furthermore, for every $M\in \bimod{T,S}$ and every $X\in\alg{T}$
there exist maps
$$\widetilde{\alpha}\colon M\circ_T EX \ra E(M\circ_S X)$$
which are natural in $M$ and $X$, and coherent with respect to
compositions of functors.
\end{cor}

\begin{exam}\label{exam-prod-pres-func}\forcepar
\begin{enumerate}
\item [(1)] Let $K$ be a fixed simplicial set, and let $X^K$ mapping
complex from $K$ to $X$.
Then $X\mapsto X^K$ commutes with  products and there are natural maps
$X\ra X^K$ induced by the projection $K\ra 1$.  Then
\eqref{prop-prod-pres-map} says that our functor $F\colon
\sSet^\indset{I}\ra \sSet^\indset{J}$ is a \emph{simplicial functor}, and
\eqref{cor-prod-pres-maps-on-alg} says that
$X^K$ is a $T$-algebra if $X$ is.

\item [(2)]  Define $E(X) = \Sing\lvert X\rvert$, the singular complex
of the geometric realization of the underlying simplicial sets;
this commutes with products and admits a natural
transformation $\Id\ra E$.  Then \eqref{cor-prod-pres-maps-on-alg} says that
$E$ lifts to all categories of $T$-algebras.

\item [(3)] Similarly, if $E(X) = \Ex^\infty(X)$, the
functor of \cite{kan;c-s-s-complexes}, then
\eqref{cor-prod-pres-maps-on-alg} says that $E$ lifts to all
categories of $T$-algebras.
\end{enumerate}
\end{exam}

\section{Simplicial algebras and s-free maps}\label{sec-s-free-maps}

In this section we explain the notion of an \emph{$s$-free
map;} this terminology 
is due to Goerss and Hopkins
\cite{goerss-hopkins-resolutions-in-model-cats}; it is essentially
what Quillen calls a \emph{free} map
\cite[II.4]{homotopical-algebra}.  

\subsection{Free degeneracy diagrams} 

Let $\simpind$ denote the category of finite totally ordered sets of
the form $[n]=\{0,\dots,n\}$ and
weakly monotone maps between them.  The category of simplicial objects
in $\cat{C}$ 
is just the category of functors $\simpind^\op\ra \cat{C}$.
The \dfn{degeneracy category} $\degencat$ is the subcategory of $\simpind$
consisting of all \emph{surjective} maps.  A \dfn{degeneracy diagram} in
$\cat{C}$ is a
functor $F\colon \degencat^\op\ra \cat{C}$. 

A degeneracy diagram $K\colon \degencat^\op\ra\Set$ is \dfn{free} if
there exist sets $L_n\subset K_n$ and an isomorphism of degeneracy
diagrams 
$$K_n \approx \coprod_{\sigma\colon[m]\ra[n]\in\degencat} L_m.$$
That is, $K$ is a left Kan extension of $L\colon \N^\op\ra \Set$ along
the inclusion 
$\N^\op\ra \simpind^\op$ sending $n\mapsto [n]$.

It is well known that if $X$ is a simplicial set, then the
underlying
degeneracy diagram of $X$ is free.  
More is true.  Let $\simpind_0\subset\simpind$ denote the subcategory
consisting of those morphisms $\delta\colon [m]\ra [n]$ such that
$\delta(0)=0$.  (A functor $\simpind_0\ra\cat{C}$ is precisely an
augmented simplicial object in $\cat{C}$ with a contracting homotopy.)
Note that $\degencat\subset\simpind_0$.
\begin{lemma}\label{lemma-free-degen}\forcepar
\begin{enumerate}
\item [(1)]
If $X\colon \simpind_0^\op\ra\Set$, then the underlying degeneracy
diagram of $X$ is free.
\item [(2)]
Suppose $Y\subset X$ is an inclusion of degeneracy diagrams of sets,
and that $X$ free.  Then $Y$
is free if and only if  
for all $x\in X_n$ and $\sigma\colon [m]\ra[n]\in\degencat$, $\sigma(x)\in 
Y_m$ implies that $x\in Y_n$.
\end{enumerate}
\end{lemma}
\begin{proof}
Let $X$ be as in (1).  For this proof, we will write simplicial
operators as acting on the right.
Say that an $x\in X_n$ in \dfn{non-degenerate}
if it is not of the form $y\sigma$ for some non-identity
$\sigma\colon [n]\ra[m]\in\degencat$ and some $y\in X_m$.  We claim:
if $x\in X_k$, $x'\in X_\ell$ are non-degenerate elements such that
$x\sigma=x'\sigma'\in X_n$ for some $\sigma, \sigma'\in \degencat$,
then
\begin{enumerate}
\item [(a)] $k=\ell$ and $x=x'$, and
\item [(b)] $\sigma = \sigma'$.
\end{enumerate}
From this claim it will follow that for each $y\in X_n$ there is a unique
non-degenerate $x\in X_m$ and a unique $\sigma\colon [n]\ra [m]\in
\degencat$ such that $y=x\sigma$; that is, 
the underlying degeneracy diagram of
$X$ is \emph{free} on the non-degenerate elements, proving (1).  To
prove the claim, observe that there exist
$\delta,\delta'\in\simpind_0$ such that $\sigma\delta=\id_{[k]}$ and
$\sigma'\delta'=\id_{[\ell]}$.  Then $x'\sigma'\delta=x\sigma\delta=x$ and
$x\sigma\delta'=x'\sigma'\delta'=x'$.  Any map in $\simpind_0$ must
factor uniquely in 
the form $\delta_1\sigma_1$ for an injective $\delta_1$ and surjective
$\sigma_1$; this fact applied to $\sigma\delta'$ and $\sigma'\delta$
together with the non-degeneracy of $x$ and $x'$ implies that
$\sigma\delta'=\sigma'\delta=\id$ and hence that $x=x'$, proving (a).  To
get (b), observe that the same argument shows that $\sigma$ and
$\sigma'$ must admit exactly the same elements of $\simpind_0$ as
right inverses, 
and it is easy to derive (b) from this.

To show (2) observe that $X$, being
free, is a disjoint union of free degeneracy diagrams on one generator (in
various degrees), and that a free degeneracy diagram on one generator
has no non-trivial free 
sub-diagrams. 
\end{proof}

\subsection{$s$-free morphisms}

We say a morphism $f\colon X\ra Y\in \alg{T}$ is \dfn{$s$-free} if,
after restricting from $\Delta$ to the degeneracy category, there is
an isomorphism 
$$Y\approx X\coprod^{\alg{T}} T(K),$$
where $K$ is a free degeneracy diagram in $\indset{I}$-graded sets.  
This
means that for each 
$n\geq0$, $Y_n \approx X_n \amalg^{T_n} T_nK_n$, and the $K_n$'s are
closed under degeneracy operations.  (The complication here is that
each level $Y_n$ in the simplicial algebra is an object in a
\emph{different} category for each $n$).  

An object $X\in \alg{T}$ is said to be \dfn{$s$-free} if the map from
the initial object to $X$ is $s$-free.  Note that $f\colon X\ra
Y\in\alg{T}$ 
is an $s$-free morphism if and only if $Y$ is an $s$-free object in
the comma category
$X\backslash\alg{T}\approx \alg{T_X}$.  

\begin{prop}\label{prop-std-res-is-s-free}
Let $X$ be a simplicial $T$-algebra, and define a simplicial object $Y$ in
$\alg{T}$ by $[n]\mapsto Y_{n,*}\approx T^{n+1}X$.  Then
$\diag(Y)\in\alg{T}$ is 
$s$-free. 
\end{prop}
\begin{proof}
We have that $Y_{n,n}\approx (T_n)^{n+1}X_n \approx T_n((T_n)^nX_n)$;
thus, we must show that $[n]\mapsto (T_n)^nX_n$ is a free degeneracy
diagram of $\indset{I}$-graded sets.  First, suppose that $T$ and $X$
are a \emph{discrete} theory and algebra.  The
degeneracy diagram $\degencat\ra \Set^\indset{I}\colon [n]\mapsto T^nX$
extends to a functor $\Delta_0\ra \Set^\indset{I}$, using the fact that
$T$ is a 
monad and and $X$ and algebra: the ``face'' maps are given by
$T^i\mu_TT^{n-i-2}\colon T^{n}X\ra T^{n-1}X$ and $T^{n-1}\psi_X\colon
T^nX\ra T^{n-1}X$.  

Since the extension from a $\degencat$-diagram to a
$\simpind_0$-diagram is natural in $T$ and $X$, we see that
$[n]\mapsto (T_n)^nX_n$ 
is the ``diagonal'' of a \emph{simplicial} object in
$\Delta_0$-diagrams, and in particular it is a $\Delta_0$-diagram, and
the result follows using \eqref{lemma-free-degen} (1). 
\end{proof}

\section{Homotopy theory of algebras}\label{sec-ho-theory-of-alg}

In this section we describe a model category structure on the category
of simplicial algebras over any $\indset{I}$-sorted theory $T$ based
on simplicial 
sets.   The model category structure we construct coincides with those
constructed in \cite{schwede-stable-ho-alg-theor} and
\cite{badzioch-alg-theories}.  

\subsection{Closed model category structure}

Let $T$ be an $\indset{I}$-sorted theory over $\sSet$.  Recall that
there is a forgetful
functor $\alg{T}\ra \sSet^\indset{I}$.  Write $U_\alpha\colon
\alg{T}\ra \sSet$ for the underlying simplicial set corresponding to
$\alpha\in \indset{I}$.

We say that a morphism $f\colon X\ra Y$ is a \dfn{strong retract} of
$g\colon X\ra Y'$ if $f$ is a retract of $g$ in the category of
objects under $X$.

\begin{thm}\label{thm-closed-model-cat}
The category $\alg{T}$
admits a simplicial  model category structure in which $f\colon
X\ra Y\in\alg{T}$ is
\begin{enumerate}
\item [(1)] a fibration or a weak equivalence if and only if each
$U_\alpha(f), \alpha\in \indset{I}$ is a fibration or weak equivalence
of simplicial sets, 
and

\item [(2)] a cofibration if and only if it is a strong retract of an
$s$-free map.
\end{enumerate}
Furthermore, this model category is right proper.
\end{thm}

Let $\phi\colon S\ra T\in s\theories(\indset{I})$ be a morphism of
$\indset{I}$-sorted simplicial theories. 
\begin{cor}\label{cor-adjoint-pair-is-quillen-pair}
The
induced adjoint pair $\phi_*\colon \alg{S}\rightleftarrows
\alg{T}\noloc \phi^*$ \eqref{prop-adjoint-pair-between-algebras} is a
Quillen pair 
between the corresponding model categories.
\end{cor}
\begin{proof}
The right adjoint $\phi^*$ is the identity on the underlying
simplicial sets, and hence preserves weak equivalences and
fibrations, and thus the left adjoint preserves cofibrations.
\end{proof}

\begin{exam}
The categories $\sSet^\indset{I}$ of graded simplicial sets admit a
model category structure in which a map is a fibration, cofibration,
or weak equivalence if it is such in each $\indset{I}$-grading.
\end{exam}
\begin{exam}\label{exam-s-theories-are-model-cat}
The category $s\theories$ of simplicial theories is a
category of algebras over an $\N$-sorted theory, and so
admits a simplicial  model category structure; similarly for
categories of bimodules over such theories.
More generally, the category $s\theories(\indset{I})$ of
$\indset{I}$-sorted simplicial theories admits a simplicial  model category
structure, as do categories of bimodules over simplicial multi-sorted theories.
\end{exam}

We will only sketch the proof of \eqref{thm-closed-model-cat}; it is
an instance of the ``small object argument'', which was already used
by Quillen \cite{homotopical-algebra} for the case of
simplicial algebras\margnote{Include note on Quillen's error.} over a
discrete theory.  (A more recent exposition of Quillen's proof for
simplicial 
algebras is \cite[II.5]{goerss-jardine-simplicial-book}.)  We note that 
the statement about right properness follows from the fact that
pullbacks, fibrations, and weak equivalences are created by the
$U_\alpha$'s, and that $\sSet$ is right proper.  The fact that
$\alg{T}$ is a \emph{simplicial} model category follows by a
straightforward argument using \eqref{cor-prod-pres-maps-on-alg} and
\eqref{exam-prod-pres-func} (1), together with the fact that graded
simplicial sets are a simplicial model category.

To apply the small object argument, we must name sets of
``generating cofibrations'' and ``generating trivial cofibrations''.
In our case we can take as generating cofibrations the set of maps
$$T(K\times \partial\Delta[n]) \ra T(K\times \Delta[n]),\qquad K\in
\ob \fsets/\indset{I},\; n\geq0,$$
and as generating trivial cofibrations the set of maps
$$T(K\times \Lambda^k[n]) \ra T(K\times \Delta[n]),\qquad
K\in\ob\fsets/\indset{I},\; n\geq k\geq0,$$
where $\Delta[n]\supset \partial\Delta[n]\supset \Lambda^k[n]$ are the
standard $n$-simplex, its boundary, and its $k$-th ``horn''.
Here we regard $\fsets/\indset{I}\subset \Set^{\indset{I}}\subset
\sSet^{\indset{I}}$ as usual, and also $\sSet\subset \sSet^\indset{I}$
by the diagonal inclusion.  (We really only need to use those $K$ whose
underlying set is a singleton.)

Using the small object argument, it is straightforward to produce
factorizations $(\text{map}) = (\text{triv.\ fib})(\text{$s$-free})$.
To get factorizations 
$(\text{map})=(\text{fib.})(\text{triv.\ cof.})$ we need the
following lemma, which ensures that the putative trivial cofibrations
produced by the small object argument are in fact such.
\begin{lemma}\label{lemma-llp-with-fibs-is-we}
Suppose $f\colon X\ra Y\in\alg{T}$ is a map which has the left lifting
property with respect to all fibrations (as defined in
\eqref{thm-closed-model-cat}).  Then $f$ is a weak
equivalence.
\end{lemma}
\begin{proof}
Let $\gamma\colon \Id\ra E$ be a natural transformation of functors
$\sSet\ra \sSet$ such that $E$ is product
preserving, $E(X)$ is a fibrant simplicial set and $\gamma_X\colon
X\ra E(X)$ is a weak equivalence for 
all $X$; we can use examples~\eqref{exam-prod-pres-func} (2) or (3).
This functor $E$ extends to $\alg{T}$ by
\eqref{cor-prod-pres-maps-on-alg}. 
Now consider 
$$\xymatrix{
{X} \ar[r]_-i \ar[d]^f
& {(EY)^{\Delta[1]} \times_{EY} EX} \ar[d]_p  \ar[r]_-\pi^-{\sim}  
& {EX}
\\
{Y} \ar[r]_j \ar@{.>}[ur]
& {EY}
}$$
where the fiber product is defined using $\ev_1\colon
(EY)^{\Delta[1]}\ra EY$ and $p$ is defined using
$\ev_0\colon (EY)^{\Delta[1]}\ra EY$.  The map $p$ is a fibration: it
can be factored 
$$(EY)^{\Delta[1]}\times_{EY}EX \ra
(EY)^{\partial\Delta[1]}\times_{EY}EX \approx EY\times EX\ra EY,$$
where both maps are fibrations since $EX$ and $EY$ are fibrant.
By hypothesis,
the dotted arrow exists.
Furthermore, $\pi$ is a trivial fibration, and hence $i$
and $j$ are weak equivalences, and we can conclude that $f$ is a weak
equivalence.  
\end{proof}

\subsection{A useful lemma}

It is convenient to give here the following generalization of
\eqref{lemma-llp-with-fibs-is-we}, which is used in
\S\ref{sec-ho-invariance-prop}. 
\begin{lemma}\label{lemma-functor-of-triv-cof-is-we}
Given the hypotheses of \eqref{lemma-llp-with-fibs-is-we}, suppose
that $F\colon \alg{T}\ra \sSet^\indset{J}$ is a degreewise functor
which commutes 
with filtered colimits and reflexive coequalizers in each degree.
Then $F(f)$ 
is a weak equivalence in $\sSet^\indset{J}$.
\end{lemma}
\begin{proof}
Consider the diagram
$$\xymatrix{
{FX} \ar[r] \ar[d]^{Ff}
& {F\left((EY)^{\Delta[1]} \times_{EY} EX\right)} \ar[d] \ar[r]
& {(EFY)^{\Delta[1]}\times_{EFY} EFX} \ar[d]
\\
{FY} \ar[r] \ar[ur]
& {FEY} \ar[r]
& {EFY.}
}$$
The left-hand side is obtained by applying $F$ to the square
used in the proof of \eqref{lemma-llp-with-fibs-is-we}.  By
\eqref{prop-char-of-filt-functors-s-algebras} the functor $F$ must be
representable by some right $T$-module, and therefore the horizontal
maps on the
right-hand side are obtained using \eqref{prop-prod-pres-map} and
\eqref{cor-prod-pres-maps-on-alg}, and the right-hand square commutes.
The top 
and bottom rows of the rectangle are weak equivalences by the same
arguments as used in the proof of \eqref{lemma-llp-with-fibs-is-we},
and hence we 
conclude that $Ff$ is a weak equivalence.
\end{proof}

\section{Homotopy invariance properties}\label{sec-ho-invariance-prop}

This section is dedicated to giving criteria for functors to preserve weak
equivalences.  As a corollary \eqref{cor-ho-inv-of-simpl-algebras} of
these results, we will see that the homotopy theory of $T$-algebras depends
only on the weak homotopy type of the simplicial
theory $T$. 

\begin{thm}\label{thm-func-we-cof-obj}
Let $T$ be an $\indset{I}$-sorted theory, and $f\colon X\ra Y$ a weak
equivalence between cofibrant 
$T$-algebras, and
let $F\colon \alg{T}\ra \sSet^\indset{J}$ be a degreewise functor
which commutes 
with filtered colimits and reflexive coequalizers (i.e., a right
$T$-module).  Then $F(f)$ 
is a weak equivalence.
\end{thm}
\begin{proof}
If $f$ is a trivial cofibration, this is
\eqref{lemma-functor-of-triv-cof-is-we}.  The theorem follows using
\eqref{lemma-ken-brown}.
\end{proof}

\begin{prop}\label{prop-we-functors-give-we}
Let $A\ra B\in{\sSet}^\fsetsdind{I}{J}$ be a weak equivalence, and
let $X\in{\sSet^\indset{I}}$.  Then the induced map $A(X)\ra B(X)$ is
a weak equivalence in $\sSet^\indset{J}$.
\end{prop}
\begin{proof}
We can first reduce to the case when $X$ is a \emph{discrete} graded
simplicial set, using the diagonal principle \eqref{subsec-notations}
and the fact that $A(X)$ (and
similarly $B(X)$) 
can be obtained as the diagonal of the simplicial object in
$\sSet^\indset{J}$ given by $[n]\mapsto A(X_n)$, where $X_n$ is the
$n$th simplicial degree of $X$.   Next note that it is enough to show that the
conclusion holds when $X$ is both
discrete and \emph{finite}, since every graded set is a filtered
colimit of its finite subsets, and $A(\mbox{--})$ and $B(\mbox{--})$
commute with such colimits.  Now we are done, since $A\ra B$ is
a weak equivalence exactly when $A(K)\ra B(K)\in \sSet^\indset{J}$ is
one for all $K\in \fsets/\indset{I}$.
\end{proof}

\begin{thm}\label{thm-ho-invariance-of-circle-over-with-fixed-alg}
Let $f\colon M\ra M'$ be a map of right $T$-modules.  The following
are equivalent.
\begin{enumerate}
\item [(1)] The map $f$ is a weak equivalence of right $T$-modules.

\item [(2)] For every $T$-algebra $X$ of the form $X=T(K)$ with $K\in
\fsets/\indset{I}\subset \sSet^\indset{I}$, the induced map
$M\circ_T X\ra M'\circ_T X$ is a weak equivalence.

\item [(3)] For every \emph{cofibrant} $T$-algebra $X$, the induced
map $M\circ_T X\ra M'\circ_T X$ is a weak equivalence.
\end{enumerate}
\end{thm}
\begin{proof}
The equivalence of (1) and (2) is immediate, since the $j$th graded
piece of $M\circ_T
T(K)$ is  $M(j,K)$.
Since for any $K\in \fsets/\indset{I}\subset \sSet^\indset{I}$, 
$X=T(K)$ is a cofibrant $T$-algebra, (3) implies (2).

To show that (1) implies (3), let $Y$ be a simplicial object in $\alg{T}$
defined by 
$[n]\mapsto Y_{n,*}=T^{n+1}X$; then $\diag(Y)\in\alg{T}$ is
$s$-free by \eqref{prop-std-res-is-s-free}, and hence is cofibrant,
and $\diag(Y)\ra X$ 
is a weak equivalence by the existence of a contracting homotopy.  Now
consider  
$$\xymatrix{
{\diag(M\circ_TY)} \ar[d]_g \ar@{=}[r]
& {M\circ_T(\diag Y)} \ar[r]^-{\sim} \ar[d]
& {M\circ_TX} \ar[d]^{f\circ_T X}
\\
{\diag(M'\circ_TY)} \ar@{=}[r]
& {M'\circ_T(\diag Y)} \ar[r]^-{\sim}
& {M'\circ_TX.}}$$
The maps marked $\sim$ are weak equivalences by
\eqref{thm-func-we-cof-obj}, so to show that $f\circ_T X$ is a weak
equivalence it suffices to show that $g$ is.  By the diagonal
principle \eqref{subsec-notations}, it suffices to show that $M\circ_T
T^{n+1}X \approx M\circ 
T^nX \ra M'\circ_T T^{n+1}X\approx M'\circ T^nX$ is a weak equivalence
for $n\geq0$; this is \eqref{prop-we-functors-give-we}. 
\end{proof}

\begin{cor}\label{cor-ho-invariance-of-circle-over-with-fixed-rt-mod}
Let $f\colon X\ra X'$ be any weak equivalence of $T$-algebras.  Then for
any \emph{cofibrant} right $T$-module $M$, the induced map $M\circ_T
X\ra M\circ_T X'$ is a weak equivalence.
\end{cor}
\begin{proof}
The functors $\mbox{--}\circ_T X, \mbox{--}\circ_T X'\colon
\rtmod{T}\ra\sSet^\indset{I}$ are represented by an
appropriate bimodules $N_X$ and $N_{X'}$, as described in
\eqref{rem-circ-on-right-has-bimod}.
We claim that the map $N_X\ra N_{X'}$
induced by $f$ is a weak equivalence, which means that we can derive
the corollary as a special case of
\eqref{thm-ho-invariance-of-circle-over-with-fixed-alg}.  To see that
$N_X\ra N_{X'}$ is a 
weak equivalence, it suffices to show that it induces a weak equivalence when
applied to a free ``algebra'', by
\eqref{thm-ho-invariance-of-circle-over-with-fixed-alg}.  Translated,
this means that we must show 
that $M\circ_T X\ra M\circ_T X'$ is a weak equivalence when $M$ is a
free right $T$-module.  In fact, this is the case whenever $M\approx
A\circ T$ for some $A\in\sSet^\fsetsdind{I}{*}$, by
\eqref{prop-we-functors-give-we}, and so is in particular true for
free objects. 
\end{proof}

\begin{rem}\label{rem-derived-circle-over}
If $M$ is a $T,S$-bimodule, then
\eqref{thm-func-we-cof-obj} 
implies
that 
the induced functor $M\circ_S\mbox{--}\colon \alg{S}\ra\alg{T}$
preserves all weak equivalences between cofibrant $S$-algebras.
Therefore, there is an induced \emph{left derived functor}
$M\circ_S^{\mathbf{L}}\mbox{--}\colon \Ho\alg{S}\ra \Ho\alg{T}$.
Similar considerations show that if $X$ is an $S$-algebra, then the
induced functor $\mbox{--}\circ_S X\colon \bimod{T,S}\ra \alg{T}$
preserves all weak equivalence between all bimodules which are
cofibrant as right $S$-modules, and hence induces a left derived
functor $\mbox{--}\circ_S^{\mathbf{L}} X\colon \Ho\bimod{T,S}\ra
\Ho\alg{T}$.  

Furthermore, \eqref{thm-ho-invariance-of-circle-over-with-fixed-alg}
and \eqref{cor-ho-invariance-of-circle-over-with-fixed-rt-mod} show
that the two ways of defining $M\circ_S^{\mathbf{L}}X$ are isomorphic
in $\Ho\alg{T}$; that is, there is a well-defined \emph{derived pairing}
$\mbox{--}\circ_S^\mathbf{L}\mbox{--}\colon
\Ho\bimod{T,S}\times\Ho\alg{S}\ra \Ho\alg{T}$.  
\end{rem}

\begin{cor}\label{cor-ho-inv-of-simpl-algebras}
Let $\phi\colon S\ra T$ be a morphism of simplicial
$\indset{I}$-sorted theories.  Then
the induced Quillen adjoint pair \eqref{cor-adjoint-pair-is-quillen-pair}
$$\phi_*\colon \alg{S}\rightleftarrows \alg{T}\noloc \phi^*$$
is a Quillen equivalence if and only if
$\phi$ is a weak equivalence 
of theories.
\end{cor}
\begin{proof}
First, note that the pair is a Quillen equivalence if and only if the
adjunction map $X\ra \phi^*\phi_*X$ is a weak equivalence for every
cofibrant $S$-algebra $X$.  This is because, given $f\colon \phi_*X\ra
Y\in \alg{T}$, the adjoint map factors
$$X \ra \phi^*\phi_*X \xra{\phi^*f} \phi^*Y,$$
and $\phi^*f$ is a weak equivalence if and only if $f$ is.
The result now follows from
\eqref{thm-ho-invariance-of-circle-over-with-fixed-alg}, since the
adjunction map is isomorphic to $S\circ_S X\ra T\circ_S X$.
\end{proof}

\section{A criterion for properness}\label{sec-proper-criterion}

In this section we give a criterion for a category of simplicial
algebras over a theory to be left proper.  The proof is adapted with 
some changes from an argument of
Dwyer and Kan 
\cite[\S8]{dwyer-kan-simplicial-localizations}, who use it to show that
simplicially enriched categories with a fixed object set form a proper
model category.  

\begin{thm}\label{thm-properness-iff-coproduct}
Let $T$ be an $\indset{I}$-sorted simplicial theory.  
The following are equivalent.
\begin{enumerate}
\item [(1)] The model category $\alg{T}$ is proper.
\item [(2)] For each finite $\indset{I}$-graded set $K\in
\fsets/\indset{I}\subset \Set^\indset{I}\subset \sSet^\indset{I}$, 
the functor $\alg{T}\ra \alg{T}$ 
given by $X\mapsto X\amalg^{T} T(K)$ carries weak
equivalences to weak equivalences.
\end{enumerate}
\end{thm}
\begin{rem}\label{rem-singleton-properness-criterion}
Note that it suffices in condition (2) of
\eqref{thm-properness-iff-coproduct} to take only those $K$ whose
underlying set is singleton.  In particular, if $\indset{I}$ is
singleton, then the theorem says that $\alg{T}$ is proper if and only
if the functor $X\mapsto X\amalg^T T(\fset{1})$ preserves weak
equivalences. 
\end{rem}

\begin{proof}
We have already seen that $\alg{T}$ is always right proper
\eqref{thm-closed-model-cat}, so we need only 
consider left properness.
That (1) implies (2) follows by observing that if $f\colon X\ra
Y\in\alg{T}$, then the square
$$\xymatrix{
{X} \ar[r] \ar[d]_f
& {X\amalg^T T(K)} \ar[d]_{g}
\\
{Y} \ar[r]
& {Y\amalg^T T(K)}
}$$
is a pushout square in $T$-algebras in which the top arrow is a
cofibration; properness implies that $g$ is a weak
equivalence if $f$ is.

To show (2) implies (1), we must show that for any cofibration
$i\colon U\ra V$,
\begin{enumerate}
\item [(*)]
the functor $\mbox{--}\amalg^T_U V\colon
U\backslash\alg{T}\ra V\backslash\alg{T}$ carries weak equivalences to
weak equivalences.  
\end{enumerate}
We proceed by a series of reductions.  First, it
suffices to show (*) when $i$ is an $s$-free map, since cofibrations are
strong retracts of such.  

Next, it suffices to show (*) for $i$ of the form $T(j)\colon
T(K)\ra 
T(L)$ where $j\colon K\ra L$ is an inclusion of $\indset{I}$-graded
simplicial sets.  This is because any $s$-free map can be written as a
directed colimit of a series of maps, each of which is a pushout along
a map of the form $T(j)$, and because weak equivalences are preserved by
directed colimits.

Define $\EuScript{B}(X,U,V)$ to be the simplicial object in $\alg{T}$
given by 
$$[n]\mapsto \EuScript{B}_n(X,U,V) = X\amalg^T\left( \coprod_n^T U
\right)\amalg^T 
V.$$
We claim that if $i=T(j)\colon T(K)\ra T(L)$, then the evident augmentation
$\diag\EuScript{B}(X,U,V)\ra X\amalg_U^T V$ is a weak equivalence.  In
fact, in each \emph{internal} degree $m$ we have that $L_m\approx
K_m\amalg K_m'$ for some $\indset{I}$-graded set $K_m'$, and thus 
$$[n]\mapsto \EuScript{B}_n(X_m,U_m,V_m) \approx X_m
\overset{T_m}{\amalg}\coprod_n^{T_m}U_m \overset{T_m}{\amalg}V_m
\approx X_m \overset{T_m}{\amalg}
T_m(\coprod_n K_m\amalg K_m\amalg K_m'),$$
which augments to $X_m\amalg^T_{T(K_m)}T(L_m)\approx X_m\amalg^{T_m}
T_m(K_m')$.  There is an evident contracting homotopy using the
inclusion $K_m'\ra K_m\amalg K_m'$, showing that
$\EuScript{B}(X_m,U_m,V_m) \ra (X\amalg^T_U V)_m$ is a weak
equivalence of 
(graded) simplicial sets,
and hence the claim follows using the diagonal principle
\eqref{subsec-notations}.

Next, it suffices to show (*) for $i$ of the form $T(\fset{0})\ra
T(K)$ for $K\in \sSet^\indset{I}$; that is, to show that the functor
$X\mapsto 
X\amalg^TT(K)$ preserves weak 
equivalences.  This follows using the diagonal principle and the above
claim, since then for $n\geq0$ each $X\mapsto X\amalg^T T(\coprod_n
K\amalg L)$ must 
preserve weak equivalences.

Next, it suffices to show (*) for $i$ of the form $T(\fset{0})\ra
T(K)$ where $K$ is a discrete graded simplicial set; this follows by
another application of 
the diagonal principle 
to $[n]\mapsto X\amalg^T T(K_n)$, the diagonal of which is $X\amalg^T
T(K)$.  

The theorem now follows using the fact that $X\amalg^T T(K)$, with $K$
discrete, is a filtered colimit over the diagram of all finite
subobjects of $K$, and that weak equivalences are preserved by filtered
colimits. 
\end{proof}

\section{Free theories and trees}\label{sec-free-theories-and-trees}

In this section we give the explicit construction of a free theory
over graded sets, and use this to derive some results needed for the
proof of \eqref{thm-cof-of-theories-cof-rt-module}.  Essentially, we
show \eqref{prop-coprod-free-theory-decomp} that a coproduct of two
\emph{free} theories is 
free as a right module over one of these theories.  That free theories
may be described in terms of trees is an observation of Boardman
\cite{boardman-language-trees},
\cite{boardman-vogt-homotopy-invariant-structures}.  
The point of view we take here is that free theories are essentially
the same as  free \emph{operads} (more precisely, free $\Sigma$-operads,
i.e., ones in which 
symmetric groups do not act), which can also be described using trees.
Our definitions of trees are based on 
those of \cite{getzler-jones}, and on ones given in an early version of
\cite{goerss-hopkins-andre-quillen-simplicial-operads}.  

\subsection{Trees}

A \dfn{totally ordered tree} $\tree{T}$ (or simply \dfn{tree}) is an oriented
contractible graph which\margnote{cite Goerss-Hopkins?  Getzler-Jones?
Boardman? where?} 
\begin{enumerate}
\item [(1)] has a (possibly empty) finite set of vertices, such that

\item [(2)] each vertex has a (possibly empty) finite totally ordered
set of input edges,

\item [(3)] each vertex has exactly one output edge, and

\item [(4)] there is exactly one edge of $\tree{T}$ which is not the output
edge of a vertex.
\end{enumerate}
Let $\inputs(v)$ denote the ordered set of input edges of a vertex $v$,
and let $\outputs(v)$ denote the unique output edge.  The external edges of
a tree $\tree{T}$ consist of a unique output edge $\outputs(\tree{T})$ and a set of
input edges $\inputs(\tree{T})$, which acquires a total ordering in an evident
way from the orderings of the $\outputs(v)$.  The output edge of a
tree is not an 
input edge, \emph{except} for the case 
of a tree which has an empty set of vertices; this is called the
\dfn{trivial tree}, and it has a unique edge. 

We fix a total ordering of each finite set $\fset{n}\in \fsets$,
so that there is a unique order preserving bijection between
$\inputs(v)$ (resp.\ $\inputs(\tree{T})$)
and some $\fset{n}$, making it convenient to identify these sets when
necessary. 

There is an evident notion of isomorphism of trees, and we will
identify isomorphic trees.

Let $\indset{I}$ be a set.  An \dfn{$\indset{I}$-tree} is a tree $\tree{T}$
together with a choice of an element $i(e)\in\indset{I}$ for each edge $e$
of $\tree{T}$; in other words, the set of edges of $\tree{T}$ is an
$\indset{I}$-graded set.  To each vertex of an $\indset{I}$-tree one
can associate an 
element $i(v)\in \N(\indset{I}) \approx \ob(\indset{I}\times
\fsets/\indset{I})$, namely the pair $(i(\outputs(v)),
i|_{\inputs(v)}\colon \inputs(v)\ra 
\indset{I})$.  Similarly, to an $\indset{I}$-tree there is an
associated element $i(\tree{T})\in \N(\indset{I})$, namely the pair
$(i(\outputs(\tree{T})), i\colon \inputs(\tree{T})\ra \indset{I})$.  

Let $A\in \Set^{\N(\indset{I})}$.  An \dfn{$A$-labelled
$\indset{I}$-tree} is a 
tree $\tree{T}$ together with a choice, for each vertex $v$ of $\tree{T}$, of an
element $a(v)\in A(i(v))$.  The set of isomorphism classes of
$A$-labelled trees is naturally 
a $\N(\indset{I})$-graded set, denoted $\freetrees A$, with the $K\in
\N(\indset{I})$ graded piece isomorphic to
$$(\freetrees A)(K) \approx \coprod_{\overset{\text{trees $\tree{T}$,}}{i(\tree{T})=K}}
\prod_{\overset{\text{vertices}}{\text{$v$ of $\tree{T}$}}} A(i(v)).$$

If $\tree{T}$ is an $A$-labelled $\indset{I}$-tree with input edges
$\inputs(\tree{T})$, and if for 
each $k\in \inputs(\tree{T})$ the
$\tree{T}_1,\dots, \tree{T}_n$ are $A$-labelled $\indset{I}$-trees such that
$i(\outputs(\tree{T}_k)) = i(k)$, then we can form a tree $\tree{T}[\tree{T}_1,\dots,\tree{T}_n]$
by \dfn{grafting} $\tree{T}_k$ at the edge $k$, obtaining a new $A$-labelled
$\indset{I}$-tree.  

\subsection{Description of free theories by trees}

Suppose $\freetheory\colon \Set^{\fsetsgr{I}} \ra
\theories(\indset{I})$ (the \dfn{free theory} functor) and
$\freeser \colon \Set^{\fsetsgr{I}} \ra \Set^{\fsetsind{I}}$ (the
\dfn{free series} functor, as in \eqref{subsec-free-series})
denote
the left adjoints to the corresponding forgetful functors.

For $A,B\in \Set^{\N(\indset{I})}$, define $A*B\in
\Set^{\N(\indset{I})}$ by
$$(A*B)(i,f\colon \fset{n}\ra\indset{I}) = 
\!\!\!\!\!\!\coprod_{g\colon \fset{m}\ra\indset{I}\in
\fsets/\indset{I}}\!\!\!\!\!\!  A(i,g) \times \left(
\coprod_{\;\underset{\text{weak monot.}}{h\colon
\fset{n}\ra\fset{m}}\;}\prod_{k\in\fset{m}} B(g(k), 
f|h^{-1}(k))\right),$$  
where the second coproduct is taken over the set of \emph{weakly monotone}
maps $h\colon \fset{n}\ra \fset{m}$ in $\fsets$ (i.e., $i\leq j$
implies $h(i)\leq h(j)$), and
$h^{-1}(k)\subset\fset{n}$ is 
identified bijectively with an object of $\fsets$ via the ordering induced
as a subset of $\fset{n}$.  Let $\delta\in\Set^{\N(\indset{I})}$ denote the
object with 
$$\delta(i, f\colon\fset{n}\ra \indset{I}) = \begin{cases} * &\text{if
$\fset{n}=\fset{1}$ and $f(1)=i$,} \\ \varnothing
&\text{otherwise.}\end{cases}$$ 
(If $\indset{I}$ is singleton, these become
$$(A*B)(n) = \coprod_m A(m) \times
\!\!\!\!\!\!\coprod_{i_1+\dots+i_m=n} \!\!\!\!\!\!
B(i_1)\times\dots \times B(i_m), \qquad \delta(n) = \begin{cases} * &
\text{if $n=1$,} \\ \varnothing & \text{otherwise.} \end{cases}.)$$
\begin{lemma}\label{lemma-monoidal-str-on-gr-sets}
\margnote{prove this?}
The category $\Set^{\N(\indset{I})}$ admits the structure of a
monoidal category, with the monoidal product given by $*$ and with
unit object $\delta$.  Furthermore, the functor $\freeser\colon
\Set^{\N(\indset{I})}\ra \Set^\fsetsind{I}$ admits the structure of a
monoidal functor for which
$I\approx \freeser \delta$ and $\freeser(A*B)\approx
\freeser A\circ\freeser B$. 
\end{lemma}
\begin{proof}
Recall from \eqref{prop-char-of-filt-functors-gr-set} that
$\Set^\fsetsind{I}$ is equivalent to a full subcategory of the
category of endofunctors on $\Set^\indset{I}$.  There is an evident
explicit 
isomorphism $\freeser A(\freeser B(X)) \approx \freeser(A*B)(X)$
natural in $X\in\Set^\indset{I}$, as can be seen by applying
\eqref{subsec-free-series}.  More explicit computations show that the
monoidal structure on $\Set^\fsetsind{I}$ restricts to
$\Set^{\N(\indset{I})}$ along $\freeser\colon \Set^{\N(\indset{I})}\ra
\Set^\fsetsind{I}$.  
\end{proof}

\begin{prop}\label{prop-free-theories-are-trees}
$\freetheory A \approx \freeser(\freetrees A)$ as objects of
$\Set^{\fsetsind{I}}$. 
\end{prop} 
\begin{rem}
The object $\freetrees A$ is nothing more than the free
$\Sigma$-operad on $A$ (cf. \cite{getzler-jones}).  Thus this proposition
relates the free $\Sigma$-operad on $A$ with the free theory on $A$. 
\end{rem}
\begin{proof}[Proof of Proposition \ref{prop-free-theories-are-trees}]
It is enough to show that $\freetrees A$ is the \emph{free monoid}
with respect to the $*$-product on $\Set^{\N(\indset{I})}$; that is,
maps $A\ra M\in \Set^{\N(\indset{I})}$ are in bijective correspondence
with maps $\freetrees A\ra M$ of monoids.  Then from
\eqref{lemma-monoidal-str-on-gr-sets} it follows formally that
$\freeser(\freetrees A)$ is the free monoid with respect to the
$\circ$-product, i.e., it is a free theory.

To make $\freetrees A$ into a monoid with respect to the $*$
structrue, let $\delta\ra 
\freetrees A$ be the map classifying the trivial trees, and let
$\freetrees A * \freetrees A\ra \freetrees A$ be the evident map
describing grafting of trees.  
\margnote{Finish
this proof.}  Now note that $\freetrees
A$ is precisely the formula for the free $\Sigma$-operad on $A$.
\end{proof}

\subsection{Essentially labelled trees}

If $\tree{T}$ is a tree, we say that $\tree{T}'\subset \tree{T}$ is a \dfn{rooted subtree}
if it is a subtree such that $\outputs(\tree{T}')=\outputs(\tree{T})$.  Given any rooted
subtree $\tree{T}'$ of $\tree{T}$ there is a unique way to write $\tree{T}$ as a graft
$\tree{T}'[\tree{T}_1,\dots,\tree{T}_n]$ for some subtrees $\tree{T}_1,\dots,\tree{T}_n$.

Let $\tree{T}\in \freetrees(A\amalg B)$.  Let $e_B(\tree{T})$ denote the minimal
rooted subtree of $\tree{T}$ which contains all of the vertices which are
labelled by $B$; if no vertices are labelled by $B$ then $e_B(\tree{T})$ is a
trivial tree.  Say that $\tree{T}\in \freetrees(A\amalg B)$ is
\dfn{$B$-essential} if $e_B(\tree{T})=\tree{T}$, and write $\esstrees(A,B)\subset
\freetrees(A\amalg B)$ for the sub-$\N(\indset{I})$-graded set of
$B$-essential trees.   We thus have shown

\begin{lemma}\label{lemma-freetrees-on-coproduct-decomp}
Every $\tree{T}\in \freetrees(A\amalg B)$ can be written uniquely as the grafting
of a $B$-essential $\indset{I}$-tree $\tree{T}'$ with $\indset{I}$-trees
$\tree{T}_1,\dots, \tree{T}_n$ labelled
only by $A$.
\end{lemma}
\begin{prop}\label{prop-coprod-free-theory-decomp}
$\freetheory(A\amalg B) \approx \freeser(\esstrees(A,B))\circ
\freetheory A$ as objects in the category of right $\freetheory
A$-modules. 
\end{prop}
\begin{proof}
Using \eqref{lemma-monoidal-str-on-gr-sets} and
\eqref{prop-free-theories-are-trees}, this amounts to showing  
that $\freetrees(A\amalg B)\approx 
\esstrees(A,B)*\freetrees A$, which is a direct translation of
\eqref{lemma-freetrees-on-coproduct-decomp}. 
\end{proof}

% Given $X\in \Set^{\indset{I}}$, let $\delta K\in \Set^{\N(\indset{I})}$
% be the object defined by 
% $$(\delta K)(i, f\colon \fset{n}\ra \indset{I}) = \begin{cases}
% K(i) & \text{if $\fset{n}=1$ and $f(1)=i$}, \\
% \varnothing & \text{otherwise.} \end{cases}$$
% In particular, if we take $X=*$ the terminal object, then $\freeser
% \delta*\in \Set^{\fsetsind{I}}$ corresponds to the 
% identity functor on $\Set^\indset{I}$.
% There is a natural map $\delta*\ra \esstrees(A,B)$, corresponding to
% the trivial trees.
\begin{prop}\label{prop-equalizer-of-free-tree-decomps}
The diagram $\esstrees(A,\varnothing)\ra\esstrees(A,B) \rightrightarrows
\esstrees(A,B\amalg B)$ is an equalizer of $\N(\indset{I})$-graded
sets, where the 
parallel maps are 
those induced by the two inclusions of $B$ into $B\amalg B$.
\end{prop}
\begin{proof}
If $\tree{T}\in \esstrees(A,B)$ has the same image under the two maps, then
it can have no vertices labelled by $B$, and hence must be a trivial
tree.  There is exactly one trivial tree for each element of $\indset{I}$,
and $\esstrees(A,\varnothing)$ contains only these.
\end{proof}

\section{Cofibrations of theories and properness}\label{sec-cofs-of-theories}

In this section we show \eqref{cor-cof-theory-proper-model-cat} that
cofibrant theories give rise to proper model categories.

\begin{thm}\label{thm-cof-of-theories-cof-rt-module}
Let $\phi\colon T\ra U$ be a cofibration between cofibrant simplicial
theories. 
Then $U$ is cofibrant as a right $T$-module. 
\end{thm}
Taking \eqref{thm-cof-of-theories-cof-rt-module} together with
\eqref{cor-ho-invariance-of-circle-over-with-fixed-rt-mod} immediately
gives
\begin{cor}\label{cor-cof-of-theories-pres-we}
If $\phi\colon T\ra U$ is a cofibration between cofibrant simplical
theories, then $\phi_*\colon \alg{T}\ra \alg{U}$ preserves all weak
equivalences. 
\end{cor}
\begin{proof}[Proof of \eqref{thm-cof-of-theories-cof-rt-module}]
We first show that it suffices to assume that $T$ is an $s$-free
theory and that $\phi$ is an $s$-free map of theories.  In fact, using
the model category structure we see that
$\phi$ is a retract of a map $\phi'\colon T'\ra U'$, where $T'$ and
$\phi'$ are $s$-free.  Then there are maps $U\ra U'\circ_{T'}T\ra U$
of right $T$-modules, and the composite of these maps is the
identity, making 
$U$ a retract of $U'\circ_{T'}T$ as a right $T$-module.  If $U'$ is
cofibrant as a right 
$T'$-module, then $U'\circ_{T'}T$ is cofibrant as a right $T$-module
(since the functor $\mbox{--}\circ_{T'}T\colon\rtmod{T'}\ra\rtmod{T}$
is the left adjoint of a Quillen pair), and hence $U$ is too.

Now suppose that $T$ and $\phi$ are $s$-free.  Thus
$T_n\approx\freetheory A_n$ and $U_n\approx\freetheory (A_n\amalg
B_n)$, where $A_n$ and $B_n$ are free degeneracy diagrams in
$\sSet^{\N(\indset{I})}$.  Then by
\eqref{prop-coprod-free-theory-decomp} we 
have $U_n\approx
\freeser(\esstrees(A_n,B_n))\circ T_n$.  Thus, it suffices to show
that $[n]\mapsto \esstrees(A_n,B_n)$ is a free degeneracy diagram in
$\Set^{\N(\indset{I})}$.  Now $\esstrees(A,B)\subset
\freetrees(A\amalg B)$, and $\freetrees(A\amalg B)$ is free by the
hypotheses that $T$ and $\phi$ be $s$-free.
By \eqref{lemma-free-degen} it suffices
to show that $\esstrees(A, B)$ is closed 
inside of $\freetrees(A\amalg B)$.
That is, if $\tree{T}\in \freetrees(A\amalg B)$ and $\sigma\in\degencat$
such that $\tree{T}\sigma \in \esstrees(A,B)$, then
$\tree{T}\in\esstrees(A,B)$.  The operator $\sigma$ acts on $\tree{T}$ by relabeling
the vertices 
of $\tree{T}$ according to the way $\sigma$ acts on $A$ and $B$ separately,
and it does not 
change the underlying shape of the tree or whether a given vertex is
labelled by $A$ or $B$; hence, if $\sigma(\tree{T})$ is
$B$-essential, then so is $\tree{T}$.
\end{proof}

Given $K\in
\Set^\indset{I}$, let $\epsilon 
K\in \Set^{\N(\indset{I})}$ denote the object with $(\epsilon
K)(\alpha, \fset{0}\ra \indset{I}) = K_\alpha$, and $(\epsilon
K)(\alpha, 
\fset{n}\ra \indset{I}) = \varnothing$ for $n>0$.  
\begin{lemma}\label{lemma-undercat-of-free-alg}
The theory $T_{T(K)}$ \eqref{subsec-undercats} is isomorphic to
$T\amalg^\freetheory 
\freetheory(\epsilon K)$, where $\freetheory(\epsilon K)$ is the free
$\indset{I}$-sorted theory on $\epsilon K$, and the coproduct is taken
in the category of $\indset{I}$-sorted
theories. 
\end{lemma}
\begin{proof}
Using the endomorphism theory technology of \eqref{subsec-undercats},
it is easy to see that $\alg{\freetheory(\epsilon K)} \approx
K\backslash\sSet^\indset{I}$.  By
\eqref{prop-alg-over-coprod-of-theories} we see that algebras over
$T\coprod^\freetheory \freetheory(\epsilon K)$ are the same as
$T$-algebras $X$ equipped with a map $K\ra X$ of graded sets, or
equivalently, the same as $T$-algebras $X$ equipped with a map
$T(K)\ra X$ of $T$-algebras.
\end{proof}

\begin{cor}\label{cor-cof-theory-proper-model-cat}
If $T$ is a cofibrant simplicial theory, then $\alg{T}$ is a proper 
model category.
\end{cor}
\begin{proof}
Suppose that $K\in\fsets/\indset{I}\subset \sSet^\indset{I}$.
By \eqref{lemma-undercat-of-free-alg}, $T\ra T_{T(K)}$ is a
cofibration between cofibrant theories, and thus 
$T_{T(K)}\circ_T\mbox{--}\colon \alg{T}\ra \alg{T_{T(K)}}\approx
T(K)\backslash \alg{T}$ carries weak equivalences to weak equivalences
by \eqref{cor-cof-of-theories-pres-we}.  Since
there is an 
isomorphism $T_{T(K)}(X) \approx X\amalg^T T(K)$ of underlying
$T$-algebras, it follows that $\alg{T}$ is proper by
\eqref{thm-properness-iff-coproduct}.  
\end{proof}

%
%   You downloaded the .tex file just to read the comments?
%   Shame on you!
%

\section{Proofs of the theorems}\label{sec-proofs-of-theorems}

\begin{proof}[Proof of Theorems A and B]
Given a simplicial theory $T$, one can construct a weak equivalence
$S\ra T$ from 
a cofibrant theory $S$, since simplicial theories are a model category
\eqref{exam-s-theories-are-model-cat}.  Then $\alg{S}$ is a proper
simplicial 
model category by 
\eqref{cor-cof-theory-proper-model-cat}, and the induced Quillen pair
$\alg{S}\rightleftarrows \alg{T}$ is a Quillen equivalence by
\eqref{cor-ho-inv-of-simpl-algebras}. 
\end{proof}

\begin{proof}[Proof of Theorem C]
Recall that $\alg{T}$ being pointed means that the initial object
$T(\fset{0})$ is
isomorphic to the terminal object, denoted $*$.  Choose $\phi\colon
S\ra T$ as in 
the proof of Theorem B, so that $\alg{S}$ is proper and is Quillen
equivalent to $\alg{T}$ via $\phi$.  The initial object in $\alg{S}$
is $S(\fset{0})$, which is not in general the terminal object.  But
since $S\ra T$ is a weak equivalence, $S(\fset{0})$ is weakly equivalent to
$T(\fset{0})\approx*$.  

Let $S_*$ denote the theory of $S$-algebras under $*$ as in \eqref{subsec-undercats}, so that
$\alg{S_*}\approx *\backslash\alg{S}$.  We have restriction functors
$\alg{T}\ra \alg{S_*}\ra\alg{S}$ factoring $\phi^*$ and hence maps
$$S\xra{\psi'} S_*\xra{\psi''} T$$
of theories factoring the weak equivalence $\phi$.  Since $\alg{S}$ is
proper, the Quillen 
pair induced by $\psi'$ is a Quillen equivalence by
\eqref{prop-undercat-char-of-proper}, and hence a weak equivalence by
\eqref{cor-ho-inv-of-simpl-algebras}.  Hence $\psi''$ is a weak
equivalence and so induces a Quillen equivalence between
$\alg{S_*}$ and $\alg{T}$.  The theorem is now proved, since
$\alg{S_*}$ is a pointed category, and is proper by
(\ref{rem-properness-conditions} (ii)). 
\end{proof}

An effective monomorphism $X\ra Y$ in a category with pushouts is a
map such that $X\ra 
Y\rightrightarrows Y\amalg_X Y$ is an equalizer.
\begin{lemma}\label{lemma-cof-theory-reg-mono}
If $T$ is a cofibrant simplicial theory, then cofibrations in
$\alg{T}$ are effective monomorphisms.
\end{lemma}
\begin{proof}
We first show that it suffices to assume that $T$ is $s$-free.  In
general, $T$ is a retract of some $s$-free $T'$.  Let $i\colon X\ra Y$
be a cofibration of $T$-algebras.  Write $X'=T'\circ_T X$ and
$Y'=T'\circ_T Y$. 
Then the diagram $X\ra
Y\rightrightarrows Y\amalg_X^TY$ is a retract of the diagram obtained
by applying $T'\circ_T\mbox{--}$ to it, which is
$X'\ra Y' \rightrightarrows Y'\amalg_{X'}^{T'}Y'$, and the map
$i'=T'\circ_T i\colon X'\ra Y'$ is a cofibration of $T'$-algebras.  If we 
know that $i'$ is an effective monomorphism, then this diagram is an
equalizer, and so is any retract of it, whence $i$ is an effective
monomorphism.

Now assume $T$ is $s$-free.  We can also assume that $i$ is an
$s$-free map, since retracts of effective monomorphisms are again
effective monomorphisms.  To show that $i$ is an effective mono, it
suffices to check it in each simplicial degree.  Thus, we must show
that for $A\in\Set^{\N(\indset{I})}$, $X\in \alg{\freetheory A}$, and
$K\in \Set^\indset{I}$, the
diagram 
$$X \ra X\amalg^{\freetheory A} (\freetheory A)(K) \rightrightarrows X
\amalg^{\freetheory A}(\freetheory A)(K\amalg K)$$
is an equalizer.  
Using \eqref{lemma-undercat-of-free-alg} and
\eqref{prop-coprod-free-theory-decomp} this is the same as
$$\freeser(\esstrees(A,\varnothing))\circ X \ra
\freeser(\esstrees(A,\epsilon K))\circ X \rightrightarrows
\freeser(\esstrees(A,\epsilon K\amalg \epsilon K))\circ X,$$
where $\epsilon K$ is as defined in \S\ref{sec-cofs-of-theories}, 
and the lemma now follows easily using
\eqref{prop-equalizer-of-free-tree-decomps}. 
\end{proof}

We note that the conclusion of \eqref{lemma-cof-theory-reg-mono} does
not hold for a general theory.  For a 
counterexample, take $\indset{I}$ singleton, and let $T$ be the unique
theory with $T(\fset{0})=\varnothing$ and $T(\fset{n})=*$ for $n>0$.
The category of $T$-algebras has exactly two objects: $\varnothing$
and $*$.  The unique map $\varnothing\ra *$ is a monomorphism, but is
not effective!

% On the other hand, if $T$ is a simplicial $\indset{I}$-theory for
% which the initial $T$-algebra, and hence every $T$-algebra, has a
% non-empty simplicial set in each $\indset{I}$-grading, then
% cofibrations of 
% $T$-algebras are effective monomorphisms.  To show this, first reduce
% the problem (as in the proof of the lemma) to 
% showing that for such a theory $T$ over sets, every map $X\ra
% X\amalg^T T(K)$ is an effective monomorphism.  Then by the non-emptyness
% hypothesis, there exists a map $r\colon K\ra X$ of graded sets, and this map
% can be used to show that the desired diagram is in fact a \emph{split}
% equalizer, with splitting induced by $r$.

\begin{proof}[Proof of Theorem D]
A model category $\cat{M}$ is \dfn{cellular} in the sense of
Hirschhorn \cite{hirschhorn} if it is a cofibrantly generated model
category with sets 
$I$ and $J$ of generating cofibrations and trivial cofibrations with
the property that
\begin{enumerate}
\item [(1)] the domains and codomains of the elements of $I$ are
``compact'',

\item [(2)] the domains of the elements of $J$ are ``small relative to
$I$'', and

\item [(3)] the cofibrations are effective monomorphisms.
\end{enumerate}
Axioms (1) and (2) say that mapping out of the domains and codomains of the
generators commutes with certain kinds of directed colimits (for the
precise notions, refer to \cite{hirschhorn}).  They certainly hold for
categories of algebras over a simplical theory, since in that case the
domains and codomains of the generators are ``small'' in the sense
that mapping out of them commutes with arbitrary filtered colimits.
Axiom (3) holds for a cofibrant theory by
\eqref{lemma-cof-theory-reg-mono}, giving the result for the
hypotheses of Theorem B.  If Axiom (3) holds in a model category, it
also holds in all undercategories, and this gives the result for the
hypotheses of Theorem C.
\end{proof}

%%% bibliography
%\bibliography{mybib}
%\bibliographystyle{amsalpha}
\newcommand{\noopsort}[1]{} \newcommand{\printfirst}[2]{#1}
  \newcommand{\singleletter}[1]{#1} \newcommand{\switchargs}[2]{#2#1}
\providecommand{\bysame}{\leavevmode\hbox to3em{\hrulefill}\thinspace}

\end{document}